\documentclass{article}
\title{On the Central Charge of a Factorizable\\
Hopf Algebra}
\author{Yorck Sommerh\"auser \qquad Yongchang Zhu}
\date{}

\usepackage{amsmath}
\usepackage{theorem}
\usepackage{amssymb}
\usepackage{latexsym}
\usepackage{diagramb}
\usepackage{tabularx}

\makeatletter
\renewcommand{\subsection}{\@startsection{subsection}{2}{0em}%
{\baselineskip}{-0em}{\bfseries\normalsize}}
\makeatother

\theoremstyle{plain}
\theorembodyfont{\rmfamily}

\newtheorem{thm}{Theorem}

\newtheorem{prop}[thm]{Proposition}
\newtheorem{lemma}[thm]{Lemma}

\newtheorem{corollary}[thm]{Corollary}

\newtheorem{conj}[thm]{Conjecture}
\newtheorem{pf}{Proof.}

\newtheorem{defn}[thm]{Definition}

\theoremstyle{break}
\theorembodyfont{\rmfamily}
\newtheorem{propb}[thm]{Proposition}
\newtheorem{defb}[thm]{Definition}
\newtheorem{corb}[thm]{Corollary}

\newcommand{\jac}[2]{\genfrac{(}{)}{}{}{#1}{#2}}

\newcommand{\qed}{$\Box$}

\newcommand{\GL}{\operatorname{GL}}
\newcommand{\SL}{\operatorname{SL}}
\newcommand{\PGL}{\operatorname{PGL}}

\newcommand{\Pn}{\operatorname{P}}

\newcommand{\End}{\operatorname{End}}

\newcommand{\Gal}{\operatorname{Gal}}

\newcommand{\Hom}{\operatorname{Hom}}

\newcommand{\id}{\operatorname{id}}

\newcommand{\Tr}{\operatorname{tr}}

\newcommand{\op}{\scriptstyle \operatorname{op}}
\newcommand{\cop}{\scriptstyle \operatorname{cop}}

\newcommand{\ab}{\operatorname{ab}}

\def\1{{1}}
\def\2{{2}}
\def\3{{3}}
\def\4{{4}}
\def\5{{5}}
\def\6{{6}}
\def\7{{7}}
\def\8{{8}}
\def\9{{9}}

\def\o{\otimes}

\def\da{\Delta}

\def\ea{\epsilon}

\def\sa{S}

\def\A{1}

\def\ua{u}

\def\C{{\mathbb C}}
\def\N{{\mathbb N}}
\def\Z{{\mathbb Z}}
\def\Q{{\mathbb Q}}
\def\R{{\mathbb R}}

\def\F{{\cal F}}

\def\T{{\mathbf T}}
\def\V{{\mathbf S}}
\def\CC{{\mathbf C}}
\def\E{{\mathbf E}}

\def\gt{{\mathfrak t}}
\def\gv{{\mathfrak s}}

\def\gd{{\mathfrak d}}

\def\bt{{\bar{\mathfrak t}}}
\def\bv{{\bar{\mathfrak s}}}

\def\gp{{\mathfrak g}}
\def\gm{{\mathfrak g}'}

\def\Gp{{\mathfrak G}}
\def\Gm{{\mathfrak G}'}

\def\fs{{\mathfrak f}}

\begin{document}

\maketitle

\begin{abstract}
\hspace{-5mm}For a semisimple factorizable Hopf algebra over a field of characteristic zero, we show that the value that an integral takes on the inverse Drinfel'd element differs from the value that it takes on the Drinfel'd element itself at most by a fourth root of unity. This can be reformulated by saying that the central charge of the Hopf algebra is an integer. If the dimension of the Hopf algebra is odd, we show that these two values differ at most by a sign, which can be reformulated by saying that the central charge is even. We give a precise condition on the dimension that determines whether the plus sign or the minus sign occurs. To formulate our results, we use the language of modular data.
\end{abstract}
\thispagestyle{empty}

\newpage

\tableofcontents

\newpage

\section*{Introduction} \label{Sec:Introd}
\addcontentsline{toc}{section}{Introduction}
In the first book of his ``Vorlesungen \"uber Zahlentheorie,'' which has been translated into English (cf.~\cite{Lan}), E.~Landau discusses various proofs of a result by C.~F.~Gau{\ss} (cf.~\cite{GaussSum}) on the sign of a certain sum of $n$~terms that we will later, in Paragraph~\ref{ClassGauss}, denote by~$\Gp$. In particular, he discusses a proof by I.~Schur (cf.~\cite{Schur}) and a proof by F.~Mertens (cf.~\cite{Mert}). As Landau points out, and this is confirmed by
the original sources, the starting point of Schur's proof is, in the notation of Paragraph~\ref{ClassGauss}, an equation of the form 
$\Gp \Gm = n$, where~$\Gm$ is the sum of the reciprocal values, whereas the starting point of Mertens' proof is an equation of the form $\Gp^2 = \pm n$,
or alternatively, by comparison with the first equation, $\Gm = \pm \Gp$.

Every semisimple factorizable Hopf algebra, or more generally every modular category, gives rise to numbers~$\gp$ and~$\gm$ that are in many ways analogous to the Gaussian sum~$\Gp$ and the reciprocal sum~$\Gm$. In fact, as we will discuss in Paragraph~\ref{ExamRadf}, both sets of numbers become equal for a suitable choice of the factorizable Hopf algebra, so that the results for Hopf algebras can be viewed as direct generalizations of the results for Gaussian sums. Therefore, looking at the above equations, the question arises whether we also have the equations~$\gp \gm = n$ 
and~$\gm = \pm \gp$. Proofs of the first of these two equations can be found in many places in the literature; we mention here only \cite{BakKir}, Cor.~3.1.11, p.~52 and \cite{Tur}, Exerc.~II.5.6, p.~116. For the second equation, less is known: A result known as Vafa's theorem implies that~$\gp/\gm$ is a root of unity (cf.~\cite{BakKir}, Thm.~3.1.19, p.~57). It is customary to write this quotient as an exponential of a parameter~$c$ called the central charge, using certain normalization conventions that we will explain in Paragraph~\ref{CentCharge}. Working through these conventions, the fact that~$\gp/\gm$ is a root of unity translates into the condition that the central charge is a rational number, this fact being at least one of the reasons for the appearance of the word `rational' in the term `rational conformal field theory'.

However, the question that arose from the classical case above was whether we have~$\gm = \pm \gp$, which translates into the condition that the central charge is an even integer. The first main result of this article is
that this is indeed the case if the dimension of the Hopf algebra is odd.
We conjecture that the result still holds if the dimension of the Hopf algebra is even. But what we can prove, and this is the second main result of this article, is that we have $\gp^4 = \gm^4$, which means for the central charge that it is an integer.

Let us be more specific. From the R-matrix of a semisimple factorizable
Hopf algebra over a field of characteristic zero, we can derive a special element, the so-called Drinfel'd element~$\ua$. The analogue of the Gaussian sum~$\Gp$ that we are considering is the 
number~$\gp = \chi_R(\ua^{-1})$, where~$\chi_R$ denotes the character of the regular representation. Note that in this situation~$\chi_R$ can also be described as the unique two-sided integral that satisfies~$\chi_R(\A)=n$,
where~$n$ denotes the dimension of the Hopf algebra. Similarly, the reciprocal Gaussian sum has the form $\gm=\chi_R(\ua)$. In the case where~$n$ is odd, our first main result is that~$\gm = \gp$ if~$n \equiv 1 \pmod{4}$, and~$\gm = -\gp$ if~$n \equiv 3 \pmod{4}$, so that~$\gp^2=\gm^2$.
We conjecture that this is still correct if~$n$ is even, but we can only
prove that~$\gp^4=\gm^4$ in this case, which is our second main result,
both of which together form exactly the two parts of Theorem~\ref{MainResult}.

All these results can be stated and proved entirely within the
Hopf algebra framework. However, we have chosen not to do this, but rather to take an axiomatic approach via modular data. This has at least three advantages: It is considerably more elementary, it makes clearer which properties are actually used to establish the results, and it facilitates the transfer of the results from the theory of Hopf algebras to related fields.

A modular datum essentially consists of two matrices that satisfy certain relations. Although these matrices and the relations between them pervade the entire literature on conformal field theory, it appears that the first systematic attempt to cast these properties into a set of simple axioms has been undertaken by T.~Gannon (cf.~\cite{GanModDat}, \cite{GanMoonMonst}). Our definition in Paragraph~\ref{DefModData} is a minor variant of his. In particular, we have designed our definition in such a way that the two essential matrices can be rescaled rather arbitrarily, and we have given names to the matrices: We speak of the Verlinde matrix~$\V$ and of the Dehn matrix~$\T$.
We have done this because the symbols used for them vary quite a bit over the literature, so that it appears inappropriate to speak of the S-matrix
and the T-matrix, in particular as this usage is exactly opposite to the 
original conventions, as used for example by R.~Fricke and F.~Klein (cf.~\cite{KleinFricke}, \S~II.2.2, p.~209).

Instead of using natural numbers $1,2,3,\ldots$ to index the rows and the columns of our matrices, we have used an abstract finite index set~$I$. The reason for this is that the basic operation for modular data is to take their Kronecker product, and the corresponding new index set, the Cartesian product of the two old ones, is not canonically modeled on a string of natural numbers. This operation is important because it forms the basis of a generalized Witt group for modular data, and, as in the classical case, the Gaussian sum is a homomorphism with respect to this group structure. However, we do not elaborate on this aspect in the sequel, at least not in this article. 

Properties of modular data that break the scalability are given special names: We call a modular datum normalized if the elements~$s_{oo}$ and~$t_o$ in the left upper corner of the Verlinde and the Dehn matrix are equal to~$1$. If the elements in the first column of the Verlinde matrix are positive integers, we call the modular datum integral. It is mainly this integrality property that distinguishes modular data coming from semisimple factorizable Hopf algebras from other modular data.

It is built into the axioms that every modular datum gives rise to a projective representation of the homogeneous modular group. By definition, rescaling the Verlinde and the Dehn matrix does not change this projective representation. Now every projective representation of the modular group can be lifted to an ordinary linear representation; as we discuss in Paragraph~\ref{CentCharge}, this is a direct consequence of the defining relations. To do this, however, we need to fix the scaling of the Verlinde and the Dehn matrix, which means that we have to choose two parameters, called the generalized rank~$D$ and the multiplicative central charge~$\ell$, a close relative of the additive central charge~$c$ considered above. A modular datum for which these two choices have been made is called an extended modular datum. As we illustrate in Paragraph~\ref{Semion} by an instructive example, the choice of these two parameters greatly affects the properties of the extended modular datum.
After the two parameters have been chosen, they can be used to scale the Verlinde and the Dehn matrix correctly, so that we get an ordinary 
representation of the modular group; these correctly scaled matrices are called the homogeneous Verlinde matrix~$\V'$ and the homogeneous Dehn matrix~$\T'$. In this context, let us note that it is customary in the literature to require that the square of the parameter~$D$ is equal to~$n$, the global dimension of the modular datum; however, the only property that is necessary for the lifting is the condition~$D^4=n^2$. For this reason,
we distinguish in Paragraph~\ref{CentCharge} between a rank and a generalized rank.

These general aspects of modular data that we have just discussed form the
contents of the preliminary Section~\ref{Sec:Prelim}. We do not claim any
originality here; we have just rearranged known material in a way that suits our later needs. The first section that contains new results is Section~\ref{Sec:OddExp}. We begin by introducing three notions:
We say what a congruence datum and a projective congruence datum is,
and define a Galois datum to be an integral modular datum whose entries~$t_i$ of the Dehn matrix are compatible with a certain action of a
Galois group described in Paragraph~\ref{GaloisAct} via the compatibility 
condition~$t_{\sigma.i} = \sigma^2(t_i)$. We then concentrate on the case where the exponent~$N$ is odd, and reach after several auxiliary steps our
first main result, namely Theorem~\ref{OddExp}: If the exponent of a normalized integral modular datum is odd, we have~$\gm = \pm \gp$.
If the global dimension~$n$ of the modular datum is also odd, then we can sharpen this assertion and obtain, by comparing with the classical Gaussian sum~$\Gp$, the result that~$\gm = \gp$ if~$n \equiv 1 \pmod{4}$ 
and~$\gm = - \gp$ if~$n \equiv 3 \pmod{4}$, which is stated in Theorem~\ref{ClassGauss}. Recall that for Hopf algebras, which are still
our main concern, it is a special case of Cauchy's theorem (cf.~\cite{YYY2},
Thm.~3.4, p.~26) that~$n$ is odd if and only if~$N$ is odd.
If we extend the modular datum using a rank, then Theorem~\ref{ClassGauss} means for the corresponding additive central charge~$c$
that it is an even integer satisfying~$c \equiv 0 \pmod{4}$ if~$n \equiv 1 \pmod{4}$ and~$c \equiv 2 \pmod{4}$ if~$n \equiv 3 \pmod{4}$.

In Section~\ref{Sec:EvenExp}, we turn to the case where~$N$ is even. Our main tool here is a result by~F.~R.~Beyl, which we quote in Paragraph~\ref{GroupCohom} and which asserts that the Schur multiplicator
of the reduced modular group~$\SL(2,\Z_N)$ is either trivial or isomorphic to~$\Z_2$. As we explain, it follows from this that, in contrast to the case of the unreduced modular group considered above, projective representations of the reduced modular group cannot always be completely lifted to an ordinary representation, but can be lifted up to a sign. 
In fact, we will see a nice explicit example of this phenomenon in Paragraph~\ref{Semion}. From Beyl's result and a couple of auxiliary steps,
we deduce in Theorem~\ref{4root} that, for a normalized extended projective congruence datum that is also Galois, we have $\gp^4 = \gm^4$. For the multiplicative central charge, this means that~$\ell^{24}=1$, which in turn means that the additive central charge~$c$ is an integer. As the example just mentioned also shows, the conclusion of this theorem cannot be sharpened to the assertion that~$\gp^2 = \gm^2$, as in the case of odd exponents. However, in the case of Hopf algebras, we conjecture that it can be sharpened in this way.

It is very important to emphasize that the same result, the integrality of~$c$, has been already derived before by A.~Coste and T.~Gannon from assumptions that look, and are, very similar (cf.~\cite{CostGann}, \S~2.4, Prop.~3.b, p.~9). Their argument is short and elegant, and we also use it as a part of our proof, as the reader will confirm when looking at Lemma~\ref{12root}. However, in the preceding terminology, Coste and Gannon
impose the Galois condition on the homogeneous Dehn matrix~$\T'$, whereas we, having the applications to Hopf algebras in mind, impose it on the original Dehn matrix~$\T$. To deduce their assumption from ours is therefore exactly equivalent to proving our result, as we explain in greater detail in Paragraph~\ref{4root}.

Section~\ref{Sec:EvenExp} ends with some results on the relation of the prime divisors of~$N$ and the prime divisors of~$n$. This is, of course,
inspired by Cauchy's theorem for Hopf algebras mentioned above. In Corollary~\ref{Cauchy}, we prove in particular that, for a projective congruence datum that is also Galois, we have $N \equiv 0 \pmod{4}$ if 
$n \equiv 2 \pmod{4}$. We do not think, however, that this result is the  best possible one, and close the section with a conjecture about a potential improvement.

Finally, we apply in Section~\ref{Sec:HopfAlg} all the machinery developed before to Hopf algebras. What we consider are semisimple factorizable Hopf algebras over fields of characteristic zero, which lead, as we explain, via the Drinfel'd element~$\ua$ to modular categories and integral modular data. The application of the preceding results is now rather straightforward. If the dimension~$n$ of the Hopf algebra is odd, Cauchy's theorem and Theorem~\ref{ClassGauss} now yield that 
$\chi_R(\ua) = \chi_R(\ua^{-1})$ if~$n \equiv 1 \pmod{4}$ and~$\chi_R(\ua) = - \chi_R(\ua^{-1})$ if~$n \equiv 3 \pmod{4}$. On the other hand, if~$n$ is even, Theorem~\ref{4root} yields that~$\chi_R(\ua)^4 = \chi_R(\ua^{-1})^4$. Of course, we have to verify the rather strong assumptions of this theorem, but here we can rely
on earlier work: That we are working with a projective congruence datum
follows from the projective congruence subgroup theorem for Hopf algebras 
(cf.~\cite{SoZhu}, Thm.~9.4, p.~94), and it is also known that our modular datum is Galois (cf.~\cite{SoZhu}, Lem.~12.2, p.~115). All these
results are obtained in Paragraph~\ref{MainResult}, where we also state
explicitly our conjecture that the 
equation~$\chi_R(\ua)^2 = \chi_R(\ua^{-1})^2$ also holds in the case of an even-dimensional Hopf algebra. We also expect that the dimension of such an even-dimensional semisimple factorizable Hopf algebra is always divisible by~$4$.

At the end of Section~\ref{Sec:HopfAlg}, we discuss an example that was constructed originally by D.~E.\ Radford (cf.~\cite{RadfAntiQuasi}, Sec.~3, p.~10; \cite{RadfKnotInv}, Sec.~2.1, p.~219). In this example, the group ring of a cyclic group of odd order~$n$ is endowed with a nonstandard \mbox{R-matrix} based on a primitive $n$-th root of unity. For this R-matrix, the Gaussian sum~$\gp$ becomes the classical Gaussian sum~$\Gp$, which shows that our preceding considerations really are generalizations of the classical case.

In the same way as Hopf algebras, quasi-Hopf algebras lead to modular data,
and these are also integral if their ribbon element is chosen correctly.
Therefore, our methods apply in just the same way, as we discuss 
in Section~\ref{Sec:QuasiHopfAlg}. However, our results in the case of a quasi-Hopf algebra are less strong, which is due to the fact that, although the projective congruence subgroup theorem has been carried over to quasi-Hopf algebras (cf.~\cite{NgSchauen4}, Thm.~8.8, p.~35), the Galois property of the corresponding modular datum has not yet been established. We are, however, optimistic that this will happen in the near future, and then our methods can be used to deduce that we have~$\gp^4 = \gm^4$ also in the quasi-Hopf algebra case. But this result can, and this is the main point of the discussion, not be improved: We construct
an explicit example of a quasi-Hopf algebra for which~$\gp^2 \neq \gm^2$.
As we mentioned above, we conjecture that it is impossible to construct
such an example with an ordinary Hopf algebra.

Both this quasi-Hopf algebra and the arising modular datum have been
considered before several times; however, we think that at least explicitly the connection between the two has not been made so far. Using a term from~\cite{RowStongWang}, we call the modular datum the semion datum; it is usually constructed using affine Kac-Moody algebras, and there
is a particularly simple case. The quasi-Hopf algebra that we are using is a dual group ring, endowed with a nontrivial associator with the help
of a 3-cocycle, a construction that has also been used frequently, also in the particularly simple case of a group of order~$2$, as we confirm by giving several explicit references. The fact that this quasi-Hopf algebra leads to the semion datum means that practically all our conjectures for Hopf algebras are false for quasi-Hopf algebras: Besides the fact 
that~$\gp^2 \neq \gm^2$ already mentioned, the semion datum is a projective congruence datum that is not a congruence datum in the sense of Definition~\ref{PropMod}. The semion datum also shows that the congruence properties explicitly depend on how the modular datum is extended by the choice of a generalized rank and a multiplicative central charge, and, as mentioned above, it provides an explicit projective representation of a reduced modular group that cannot be lifted to a linear representation, thus confirming the nontriviality of the corresponding Schur multiplicator.

Finally, let us say some words about the conventions that we use throughout the exposition. Our base field is denoted by~$K$, and while we always assume that it is of characteristic zero, we do not always assume that it is algebraically closed. The dual of a vector space~$V$ is denoted by $V^*:=\Hom_K(V,K)$. Unless stated otherwise, a module is a left module. If~$R$ is a ring, then~$R^\times$ denotes its group of units. Also, we 
use the so-called Kronecker symbol~$\delta_{ij}$, which is equal to~$1$ if $i=j$ and equal to~$0$ otherwise. The set of natural numbers is the set \mbox{$\N:=\{1,2,3,\ldots\}$}; in particular, $0$ is not a natural number. The symbol~$\Q_N$ denotes the $N$-th cyclotomic field, and not some field of $N$-adic numbers, and~$\Z_N$ denotes the set~$\Z/N\Z$ of integers modulo~$N$, and not some ring of $N$-adic integers.
The greatest common divisor of two integers~$k$ and~$l$ is denoted by~$\gcd(k,l)$ and is always chosen to be nonnegative. With respect to enumeration, we use the convention that propositions, definitions, and similar items are referenced by the paragraph in which they occur; an 
additional third digit indicates a part of the corresponding item.
For example, a reference to Proposition~1.1.1 refers to the first
assertion of the unique proposition in Paragraph~1.1. 

Part of the present work was carried out during a visit of the first author to the Hong Kong University of Science and Technology. He thanks the university, and in particular his host, for the hospitality.

\newpage
\section{Preliminaries} \label{Sec:Prelim}
\subsection[Modular data]{} \label{DefModData}
Although modular data have played a very prominent role in conformal field
theory for quite some time, it appears that the approach to derive their properties from a set of simple axioms was first proposed by T.~Gannon (cf.~\cite{GanModDat}, Def.~1, p.~214; \cite{GanMoonMonst}, Def.~6.1.6, p.~359). We use here a slight modification
of his system of axioms; the relation of the two systems of axioms will
become clear in Paragraph~\ref{CentCharge}.

\begin{defn}
Suppose that $K$ is a field of characteristic zero. A modular datum over~$K$ is a quintuple $(I,o,{}^*,\V,\T)$ consisting of a finite set~$I$,
a distinguished element $o \in I$, called the unit element, an involution
$$^*: I \rightarrow I,~i \mapsto i^*$$
a matrix $\V = (s_{ij})_{i,j \in I}$ with entries in~$K$, called the Verlinde matrix, and a diagonal matrix $\T = (t_i \delta_{ij})_{i,j \in I}$  with entries in~$K$, called the Dehn matrix.
We require these to satisfy the following axioms:
\begin{enumerate}
\item  \label{DefModData1}
$\V$ is symmetric, and $\T$ has finite order.

\item  \label{DefModData2}
For all $i \in I$, we have $t_{i^*}=t_i$ and $s_{io} \neq 0$. Moreover, we have $o^*=o$.

\item  \label{DefModData3}
There is a nonzero number $n \in K$ such that
$$\sum_{j \in I} s_{ij} s_{jk} = n\delta_{i,k^*}$$
for all $i,k  \in I$.

\item  \label{DefModData4}
The matrices $\T^{-1} \V \T^{-1}$ and $\V \T \V$ are proportional.

\item  \label{DefModData5}
The numbers
$$N_{i,j}^k := 
\frac{1}{n} \sum_{l \in I} \frac{s_{il} s_{jl} s_{k^*l}}{s_{ol}}$$
are nonnegative integers.
\end{enumerate}
A modular datum is called normalized if $s_{oo}=t_o=1$. It is called integral if the elements $s_{io}$ are positive integers.
\end{defn}

When we refer to elements of the base field~$K$ as integers, what we really mean is that this field element is
the image of an integer under the unique field isomorphism between~$\Q$
and the prime field of~$K$. In the sequel, we will frequently pass in this way from a rational number to a field element and vice versa without explicit mention.

Note that, with the exception of normality and integrality, the axioms of a modular datum have the property that they are still satisfied if we rescale the Verlinde matrix by a nonzero factor and the Dehn matrix by a root of unity. 

The various constants that appear in this definition have
special names. If \linebreak
$\E=(\delta_{ij})_{i,j \in I}$ is the unit matrix,
then the definition implies that there is a smallest number~$N$, called
the exponent of the modular datum, such that $\T^N=\E$. The diagonal
entries~$t_i$ of the Dehn matrix are then $N$-th roots of unity.
The normalized Dehn matrix~$\T/t_o$ then also has finite order, so that there is a smallest number~$N_o$, called
the normalized exponent of the modular datum, such that $(\T/t_o)^{N_o}=\E$. Clearly, $N_o$ divides~$N$.

The number~$n$ appearing in Axiom~\ref{DefModData3} is called the global dimension of the modular datum. Inserting $i=k=o$ into this axiom, we get that $n=\sum_{j \in I} s_{jo}^2$. For $i \in I$, we also call the number
$n_i := s_{io}$ the $i$-th dimension, so that we can reformulate
the last equation as $\sum_{j \in I} n_{j}^2 = n$. If we introduce the charge conjugation matrix 
$\CC:=(\delta_{i,j^*})_{i,j \in I}$, then Axiom~\ref{DefModData3}
can be written in the matrix form $\V^2 = n \CC$. Since $n \neq 0$, we see in particular that both~$\V$ and~$\T$ are invertible. Note also that the integrality condition means that the dimensions~$n_i$ are positive integers, and therefore the global dimension~$n$ is also a positive integer. In this context, let us point out that our definition of integrality deviates from the one in \cite{Cuntz}, Def.~3.1, where in particular all entries of the Verlinde matrix are required to be integers.

We can also determine the proportionality constant in Axiom~\ref{DefModData4}. The $(o,o)$-component of~$\V \T \V$ is the so-called Gaussian sum  $\gp := \sum_{i \in I} n_i^2 t_i$ of the modular datum, whereas the $(o,o)$-component of~$\T^{-1} \V \T^{-1}$
is $n_o/t_o^2$, so that Axiom~\ref{DefModData4} becomes
\begin{align*} \label{DefModData4WithConst}  
\gp \T^{-1} \V \T^{-1}= \frac{n_o}{t_o^2} \V \T \V
\end{align*}
which we will call the constant form of Axiom~\ref{DefModData4}.

We list some elementary consequences of our axioms that will be frequently used in the sequel:
\begin{propb}
\begin{enumerate}
\item \label{DefModDataStar}
We have $s_{i^*,j^*}=s_{ij}$ and $n_{j^*} = n_{j}$.

\item
The matrix~$\CC$ commutes with~$\V$ and~$\T$; i.e., we have
$\CC \V = \V \CC$ and $\CC \T = \T \CC$.

\item \label{DefModDataNo}
For all $i,j \in I$, we have $N_{oi}^j = N_{io}^j = \delta_{ij}$
and $N_{ij}^o = \delta_{i,j^*}$.

\item \label{DefModDataST}
For all $i,j \in I$, we have
$$s_{ij} = \frac{t_o}{t_i t_j} \sum_{k\in I} N_{ik}^j n_k t_k$$

\item \label{DefModDataRecGauss}
We have $\gp \gm = n n_o^2$, where
$\gm := \sum_{i \in I} \frac{n_i^2}{t_i}$
is the so-called reciprocal Gaussian sum. In particular, the Gaussian sum
and its reciprocal are nonzero.
\end{enumerate}
\end{propb}
\begin{pf}
\newcounter{num}
\begin{list}{(\arabic{num})}{\usecounter{num} \leftmargin0cm \itemindent5pt}
\item
Written in matrix form, the first part of the first assertion says that
$\CC \V \CC = \V$, or alternatively $\CC \V = \V \CC$. But this holds since $\CC=\frac{1}{n} \V^2$. Setting $i=o$ gives the second part of the first assertion. For the second assertion, we have just seen that~$\V$ commutes with~$\CC$, and the fact that~$\T$ commutes with~$\CC$ is equivalent to our axiom that $t_{i^*}=t_i$. The third assertion
is a direct consequence of Axiom~\ref{DefModData3}.

\item
For the fourth assertion, the $(o,l)$-component of the constant form of Axiom~\ref{DefModData4} yields
$$\frac{\gp n_l}{t_o t_l} = 
\frac{n_o}{t_o^2} \sum_{j \in I} n_{j} t_j s_{jl}$$
Using the definition of the numbers $N_{ij}^k$, we therefore get
\begin{align*}
\sum_{j\in I} N_{ij}^k n_j t_j = 
\frac{1}{n} \sum_{j,l\in I} \frac{s_{il} s_{jl} s_{k^*l}}{n_l} n_j t_j = 
\frac{t_o \gp}{n n_o} \sum_{l\in I} s_{il} s_{k^*l} \frac{1}{t_l}  = 
\frac{t_i t_k}{t_o} s_{ik} 
\end{align*}
where the last step uses the matrix identity
$\gp \V^{-1} \T^{-1} \V = \frac{n_o}{t_o^2} \T \V \T$.

\item
For the fifth assertion, we first note that the Gaussian sum is clearly nonzero, because the right-hand side in the constant form of Axiom~\ref{DefModData4} is invertible. 
If we invert this very equation, it becomes
$$\frac{1}{\gp} \T \V^{-1} \T = 
\frac{t_o^2}{n_o} \V^{-1} \T^{-1} \V^{-1}$$
As we have $\V^{-1} = \frac{1}{n} \V \CC$, we can write this in the form
$$\frac{1}{\gp} \T \V \CC \T = 
\frac{t_o^2}{n n_o} \V \T^{-1} \V$$
The $(o,o)$-component of this matrix equation is
\begin{align*}
\frac{t_o^2 n_o}{\gp} = 
\frac{t_o^2}{n n_o} \sum_{i \in I} s_{oi} \frac{1}{t_i} s_{io} 
= \frac{t_o^2}{n n_o} \gm
\end{align*}
which yields the fifth assertion.
\qed
\end{list}
\end{pf}

Note that variants of the fourth assertion appear in \cite{GanModDat}, Sec.~4, p.~232 and \cite{Tur}, Sec.~II.4.5, p.~108. The fifth assertion  should be compared with \cite{Tur}, Sec.~II.2.4, Eq.~(2.4.a), p.~83, and Exerc.~II.5.6, p.~116. As mentioned in the introduction, the fifth assertion gives a first relation between the Gaussian sum and its reciprocal; the main issue in this article will be to find more such relations.

\subsection[The fusion ring]{} \label{FusRing}
With every modular datum, we associate a ring~$\F$, called the fusion ring, as follows: As an abelian group, $\F$ is free on a basis 
$\{b_i \mid i \in I\}$.
The multiplication is given on basis elements by
$$b_i b_j = \sum_{k \in I} N_{ij}^k b_k$$
and then on the whole group by bilinear extension. We will
also frequently extend the scalars to other commutative rings~$R$ by
considering $\F_R := \F \o_\Z R$.

As we have $N_{ij}^k = N_{ji}^k$ by construction, $\F$ is clearly commutative, and it follows from Proposition~\ref{DefModData}.\ref{DefModDataNo} that~$b_o$ is a unit element. But this multiplication is actually also
associative, as we will show now. For this, we define for each $q \in I$ a group homomorphism
$\xi_q: \F \rightarrow K$ to the additive group of the field by requiring that the value on the basis elements be
$$\xi_q(b_i) = \frac{s_{iq}}{n_q}$$
This map is not only a group homomorphism:
\begin{lemma}
$\xi_q$ is a ring homomorphism.
\end{lemma}
\begin{pf}
By Axiom~\ref{DefModData}.\ref{DefModData3} and Axiom~\ref{DefModData}.\ref{DefModData5}, we have
\begin{align*}
\xi_q(b_i b_j) = \sum_{k \in I} N_{ij}^k \frac{s_{kq}}{n_q} &= 
\frac{1}{n} \sum_{k,l \in I} \frac{s_{il} s_{jl} s_{k^*l}}{n_l}  \frac{s_{kq}}{n_q} \\
&= \sum_{l\in I} \frac{s_{il} s_{jl} \delta_{ql}}{n_l n_q} 
= \frac{s_{iq}}{n_q} \frac{s_{jq}}{n_q} = \xi_q(b_i) \xi_q(b_j)
\end{align*}
where we have also used Proposition~\ref{DefModData}.\ref{DefModDataStar} for the third equality. The fact that $\xi_q(b_o)=1$ follows directly from the definition.
\qed
\end{pf}

Via this lemma, the associativity of~$\F$ follows essentially from 
the associativity of~$K$:
\begin{prop}
$\F$ is a commutative, associative, reduced ring with unit element~$b_o$. 
The algebra $\F_\Q$ is semisimple.
\end{prop}
\begin{pf}
If~$K^I$ denotes the algebra of functions on~$I$ with pointwise addition
and multiplication, then it follows from the preceding lemma that
the map
$$\F \rightarrow K^I,~b \mapsto (\xi_i(b))_{i \in I}$$
is a ring homomorphism. Because the Verlinde matrix is invertible,
this map is injective, and the map from~$\F_K$ to~$K^I$ given by extension of scalars is in fact an algebra isomorphism. Since~$K^I$ is associative and reduced, this also holds for~$\F$. Finally, it is a general fact
that a finite-dimensional commutative reduced algebra over a field is semisimple (cf.~\cite{FarbDennis}, Thm.~2.2, p.~60).
\qed
\end{pf}

The algebra isomorphism between~$\F_K$ and~$K^I$ appearing in the preceding proof can be used to transfer all properties from~$K^I$ to~$\F_K$; in particular, we can easily write down the primitive idempotents of~$\F_K$. For this, we introduce the element
$$b_A := \sum_{j \in I} b_j b_{j^*}$$
which relates to the primitive idempotents as follows:
\begin{corb}
\begin{enumerate}
\item \label{bA} 
For $i \in I$, we have $\xi_i(b_A) = \frac{n}{n_i^2} \neq 0$.

\item \label{IdemExpl} 
The element 
$p_i := \frac{1}{\xi_i(b_{A})} \sum_{j \in I} \xi_i(b_{j^*}) b_j$
is a primitive idempotent in~$\F_K$.

\item 
For $i,j \in I$, we have $\xi_j(p_i)=\delta_{ij}$.

\item
For $b \in \F_K$, we have $b p_i = \xi_i(b) p_i$.
\end{enumerate}
\end{corb}
\begin{pf}
As we have $\xi_i(b_j) = \frac{s_{ij}}{n_i}$, it follows from Axiom~\ref{DefModData}.\ref{DefModData3} that
$$\xi_i(b_{A}) = \sum_{j \in I} \frac{s_{ij}}{n_i} \frac{s_{ij^*}}{n_i} 
= \frac{1}{n_i^2} \sum_{j \in I} s_{ij} s_{ji^*}
= \frac{n}{n_i^2}$$
where we have also used Proposition~\ref{DefModData}.\ref{DefModDataStar}.
This proves the first assertion. For the following assertions, we should first point out that we have denoted the basis
element $b_j \o 1 \in \F_K = \F \o_\Z K$ again by~$b_j$, and also denoted the $K$-linear extension of~$\xi_i$ from~$\F_K$ to~$K$ again by~$\xi_i$.
The third assertion then holds because
\begin{align*}
\xi_j(p_i) = 
\frac{1}{\xi_i(b_{A})} \sum_{k \in I} \xi_i(b_{k^*}) \xi_j(b_{k}) 
= \frac{n_i^2}{n} \sum_{k \in I} \frac{s_{ik^*}}{n_i} \frac{s_{jk}}{n_j}
= \delta_{ij}
\end{align*}
by Axiom~\ref{DefModData}.\ref{DefModData3} and Proposition~\ref{DefModData}.\ref{DefModDataStar}.

For the second assertion, we denote the algebra isomorphism appearing in the proof of the preceding proposition by 
$f: \F_K \rightarrow K^I,~b \mapsto (\xi_j(b))_{j \in I}$. It follows from the third assertion that
$f(p_i) = (\delta_{ij})_{j \in I}$, which is a primitive idempotent in~$K^I$. Therefore, $p_i$ must be a primitive idempotent in~$\F_K$.
The fourth assertion follows similarly: We have
\begin{align*}
f(b p_i) = f(b) f(p_i) = (\xi_j(b))_{j \in I} \; (\delta_{ij})_{j \in I}
= \xi_i(b) \; (\delta_{ij})_{j \in I} = \xi_i(b)  f(p_i)
\end{align*}
and therefore $b p_i = \xi_i(b) p_i$.
\qed
\end{pf}

The special case $i=o$ of the fourth assertion in this corollary is worth noting: Since $\xi_o(b_j) = n_j/n_o$, we have 
$p_o = \frac{n_o}{n} \sum_{j \in I} n_j b_j$,
and the fourth assertion then yields that
$b_k p_o = \frac{n_k}{n_o} p_o$
for all $k \in I$. 

It must be emphasized that the material in this paragraph is not
at all new; it can be found in many places, among them \cite{ChariPress}, Def.~5.2.8, p.~155; \cite{GanModDat},
Def.~2, p.~216; and \cite{GanMoonMonst}, Def.~6.1.3, p.~355.
An analogous discussion in the Hopf algebra case can be found in \cite{SoZhu}, Par.~5.1 and Par.~5.3. It is also worth noting that
the rescaled Verlinde and Dehn matrices lead to the same fusion
ring, as only homogeneous ratios of the matrix elements enter into
the definition.

\subsection[Cyclotomic fields]{} \label{Cyclotomic}
By construction, the diagonal elements~$t_i$ of the Dehn matrix generate
the $N$-th cyclotomic field~$\Q_N \subset K$, and the ratios~$t_i/t_o$ generate its subfield~$\Q_{N_o}$. If the modular datum is
integral, then we see by multiplying and dividing the right-hand side of the formula in Proposition~\ref{DefModData}.\ref{DefModDataST} by~$t_o$ 
that the entries of the Verlinde matrix are in fact contained in~$\Q_{N_o}$. This has the following consequence:
\begin{prop}
For an integral modular datum, the algebra $\F_\Q$ is isomorphic
to a direct sum of subfields of the cyclotomic field~$\Q_{N_o}$.
\end{prop}
\begin{pf}
By the Wedderburn structure theorem (cf.~\cite{FarbDennis}, Thm.~1.11, p.~40), we can decompose~$\F_\Q$ in the form
$$\F_\Q = \bigoplus_{q=1}^r K_q$$
where~$K_1,\ldots,K_r$ are fields. As we saw in Paragraph~\ref{FusRing}, the algebra homomorphism $\F_\Q \rightarrow K^I,~b \mapsto (\xi_i(b))_{i \in I}$
is injective. Therefore, for every~$q \le r$ there is an $i \in I$ such 
that $\xi_i(K_q)$ is nonzero. But as $\xi_i(b_j) = s_{ji}/n_i \in \Q_{N_o}$, the image of~$\xi_i$ is contained in~$\Q_{N_o}$. Since field homomorphisms are always injective, $\xi_i$ restricts to an isomorphism between~$K_q$ and a subfield of~$\Q_{N_o}$.~\qed
\end{pf}

The preceding argument should be compared with 
\cite{BoerGoer}, App.~B, p.~302, where a similar result is proved using the Kronecker-Weber theorem. Our argument also shows that the algebra homomorphisms~$\xi_i$ arise from the Wedderburn decomposition by
projecting to a Wedderburn component~$K_q$ and then embedding it into~$\Q_{N_o}$: Because the product of elements in distinct Wedderburn components is zero, and~$\Q_{N_o}$ does not contain nontrivial zero divisors, every $\xi_i$~must vanish on all Wedderburn components~$K_q$ except one, and on this Wedderburn component it becomes an embedding into~$\Q_{N_o}$.

The maps~$\xi_i$ actually exhaust all algebra homomorphisms from~$\F_\Q$ to~$\Q_{N_o}$. This follows from Dedekind's theorem (cf.~\cite{Jac3}, \S~1.3, Thm.~3, p.~25), which implies that distinct algebra homomorphisms are linearly independent elements in the vector space
$\Hom_\Q(\F_\Q,\Q_{N_o})$ over~$\Q_{N_o}$, whose dimension is the cardinality of~$I$. Because the Verlinde matrix is invertible, it cannot have two proportional columns, which implies that the maps~$\xi_i$ are all distinct.

\subsection[The action of the Galois group]{} \label{GaloisAct}
For an integral modular datum, consider the Galois group $\Gal(\Q_{N_o}/\Q)$ of the cyclotomic field~$\Q_{N_o}$. We have seen in Proposition~\ref{Cyclotomic} that the Wedderburn components~$K_q$ in the Wedderburn decomposition
$$\F_\Q = \bigoplus_{q=1}^r K_q$$
are fields that can be embedded into~$\Q_{N_o}$. Because the Galois group 
$\Gal(\Q_{N_o}/\Q)$ is abelian, all subfields of~$\Q_{N_o}$ are invariant under its action, and therefore we get an action of the Galois group on~$K_q$, which for the same reason does not depend on the choice of the embedding. By summing up the actions on all the Wedderburn components, we get an action of~$\Gal(\Q_{N_o}/\Q)$ on~$\F_\Q$, which has the property that
$$\xi_i(\sigma.b) = \sigma(\xi_i(b))$$
for all $\sigma \in \Gal(\Q_{N_o}/\Q)$, all $i \in I$, and all $b \in \F_\Q$. This holds because
we saw in Paragraph~\ref{Cyclotomic} that the algebra homomorphisms~$\xi_i$
are exactly those which arise by projecting to some Wedderburn component and then embedding it into~$\Q_{N_o}$. 

Now $\sigma \circ \xi_i$ is again an algebra homomorphism from~$\F_\Q$ to~$\Q_{N_o}$,
and therefore there must be an index~$j \in I$ such that
$\sigma \circ \xi_i = \xi_{j}$. We denote this index~$j$ by
$\sigma.i$, and get in this way an action of~$\Gal(\Q_{N_o}/\Q)$
on the index set~$I$ with the property that
$$\sigma(\xi_i(b)) = \xi_{\sigma.i}(b)$$
for all $b \in \F_\Q$. The following proposition lists the basic properties
of the actions that we have introduced:
\begin{prop}
For an integral modular datum, we have
\begin{enumerate}
\item 
$\sigma(s_{ij}) = s_{\sigma.i,j} = s_{i,\sigma.j}$

\item 
$n_{i} = n_{\sigma.i}$

\item
$\sigma.b_i = b_{\sigma.i}$

\item
$\sigma.o=o$
\end{enumerate}
\end{prop}
\begin{pf}
\begin{list}{(\arabic{num})}{\usecounter{num} \leftmargin0cm \itemindent5pt}
\item
From Corollary~\ref{FusRing}, we get that
$$\frac{n}{n_{\sigma.i}^2} = \xi_{\sigma.i}(b_A) = \sigma(\xi_{i}(b_A))
= \sigma(\frac{n}{n_{i}^2}) = \frac{n}{n_{i}^2}$$
This implies $n_{\sigma.i}^2 = n_{i}^2$, which implies 
$n_{\sigma.i} = n_{i}$ because every~$n_i$ is a positive integer by assumption. This proves the second assertion.

\item
By definition, we have
$\xi_i(b_j) = \frac{s_{ij}}{n_i}$, so that the preceding equation implies
\begin{align*}
\frac{s_{\sigma.i,j}}{n_{\sigma.i}} = \xi_{\sigma.i}(b_j) 
= \sigma(\frac{s_{ij}}{n_i}) = \frac{\sigma(s_{ij})}{n_i}
\end{align*}
This implies $\sigma(s_{ij}) = s_{\sigma.i,j}$ by the preceding step.
By the symmetry of the Verlinde matrix, we then also have
$$\sigma(s_{ij}) = \sigma(s_{ji}) = s_{\sigma.j,i} = s_{i,\sigma.j}$$
which proves the first assertion.

\item
The equation $\sigma.b_i = b_{\sigma.i}$ will hold if we can show 
that $\xi_j(\sigma.b_i) = \xi_j(b_{\sigma.i})$ for all~$j \in I$. But this equation holds, as we have
\begin{align*}
\xi_j(\sigma.b_i) = \sigma(\xi_j(b_i)) = \sigma(\frac{s_{ji}}{n_j}) =
\frac{s_{j,\sigma.i}}{n_j} = \xi_j(b_{\sigma.i})
\end{align*}

\item
By construction, $\sigma$ acts as a unital algebra automorphism on~$\F_\Q$,
and in particular preserves the unit element, so that $\sigma.b_o=b_o$.
By the preceding step, this implies $\sigma.o=o$.
\qed
\end{list}
\end{pf}

Because~$N_o$ divides~$N$, $\Gal(\Q_{N_o}/\Q)$ is a factor group of~$\Gal(\Q_{N}/\Q)$, which implies that~$\Gal(\Q_{N}/\Q)$ also acts on~$I$ via pullback along the canonical map. This fact will be used
in the sequel without explicit mention, for example in Definition~\ref{PropMod}.\ref{GaloisDat}.

\subsection[Complex conjugation]{} \label{ComplConj}
The Galois group $\Gal(\Q_{N_o}/\Q)$ contains an element~$\gamma$ that maps
every \mbox{${N_o}$-th} root of unity to its inverse. It corresponds to complex conjugation under every embedding of the cyclotomic field into the complex numbers. As we show now, this automorphism~$\gamma$ acts on the index set as the involution that is part of the definition of our modular datum, which is still assumed to be integral: 
\begin{prop}
For all $i \in I$, we have $\gamma.i =i^*$.
\end{prop}
\begin{pf}
Note first that from Proposition~\ref{DefModData}.\ref{DefModDataStar}
we get
\begin{align*}
N_{i^*j^*}^{k^*} = 
\frac{1}{n} \sum_{l\in I} \frac{s_{i^*l} s_{j^*l} s_{kl}}{n_l} = 
\frac{1}{n} \sum_{l\in I} \frac{s_{i^*l^*} s_{j^*l^*} s_{kl^*}}{n_{l^*}} = 
\frac{1}{n} \sum_{l\in I} \frac{s_{il} s_{jl} s_{k^*l}}{n_{l}}
= N_{ij}^k
\end{align*}
which shows that the $\Q$-linear endomorphism of~$\F_\Q$ that maps the basis element~$b_i$ to~$b_{i^*}$ is an algebra homomorphism. We extend this map antilinearily to the complexification~$\F_\C$; i.e., for an element
$x = \sum_{i \in I} z_i b_i$, with~$z_i \in \C$, we define
$$x^* = \sum_{i \in I} \bar{z}_i b_{i^*}$$
where $\bar{z}_i$ is the complex conjugate of~$z_i$.

The basis $\{b_i \mid i \in I\}$ of~$\F_\C$ has a dual basis, and we denote the dual basis vector of~$b_o$ by~$\lambda$, so that we have
$\lambda(b_i) = \delta_{i,o}$. Using it, we define on~$\F_\C$
the sesquilinear form
$$\langle x, y \rangle := \lambda(xy^*)$$
Now it follows from Proposition~\ref{DefModData}.\ref{DefModDataNo}
that the basis $\{b_i \mid i \in I\}$ is orthonormal with respect to this sesquilinear form, which shows that the sesquilinear form is also
positive definite, i.e., a scalar product on~$\F_\C$. This shows in particular that $xx^*=0$ implies $x=0$ for all $x \in \F_\C$.

From Corollary~\ref{FusRing}.\ref{IdemExpl}, we see that the primitive idempotents $p_i \in \F_{K}$ are already contained in~$\F_{\Q_{N_o}}$. By choosing a fixed embedding of~$\Q_{N_o}$ into~$\C$, we can regard them as elements of~$\F_\C$. Because ${}^*$ is a ring automorphism, $p_i^*$ is again a primitive idempotent, and therefore 
$p_i \neq p_i^*$ implies $p_i p_i^* = 0$, which would imply $p_i=0$ by the preceding argument. This shows that
$p_i = p_i^*$ for all $i \in I$. Writing this out in the explicit form of these elements given in Corollary~\ref{FusRing}.\ref{IdemExpl},
this means that
$$\frac{1}{\xi_i(b_{A})} \sum_{j \in I} \xi_i(b_{j^*}) b_j = 
\frac{1}{\xi_i(b_{A})} \sum_{j \in I} \gamma(\xi_i(b_{j^*})) b_{j^*} $$
where we have used that $\xi_i(b_{A})=\frac{n}{n_i^2}$ is a rational number by Corollary~\ref{FusRing}.\ref{bA}. This shows that
$\xi_i(b_{j^*}) = \gamma(\xi_i(b_{j}))$. By the definition of the action of the Galois group, this means that 
$\xi_i(b_{j^*}) = \xi_i(\gamma.b_{j})$, which implies by Proposition~\ref{GaloisAct} that
$b_{j^*} = \gamma.b_{j} = b_{\gamma.j}$, as asserted.
\qed
\end{pf}

The preceding argument is adapted from \cite{NR2}, Sec.~2, Rem.~11, p.~1088, resp.~\cite{With}, Sec.~3, Prop.~3.1, p.~885; a slightly different argument can be found in~\cite{Beauv}, \S~6, Cor.~6.2, p.~88. Note that it implies in particular that $\sigma.i^*=(\sigma.i)^*$ for all 
automorphisms $\sigma \in \Gal(\Q_{N_o}/\Q)$, since this group is abelian.

\subsection[The Gaussian sum]{} \label{GaussSum}
As we said in the introduction and will see in Paragraph~\ref{ExamRadf} in explicit example, the Gaussian sum~$\gp$ is in general not equal to the reciprocal Gaussian sum~$\gm$. However, they differ only by a root of unity:
\begin{prop}
Let~$m$ be the cardinality of~$I$.
\begin{enumerate}
\item
We have $\gp^{2Nm}=\gm^{2Nm}$.

\item 
If the modular datum is integral, we have $\gp^{2N}=\gm^{2N}$. 

\item 
If the modular datum is integral and $N$ is even, we have $\gp^{N}=\gm^{N}$.
\end{enumerate}
\end{prop}
\begin{pf}
\begin{list}{(\arabic{num})}{\usecounter{num} \leftmargin0cm \itemindent5pt}
\item
By taking determinants on both sides of the equation $\V^2 = n \CC$, we get that
$$\det(\V)^2 = n^m \det(\CC) = \pm n^m$$
because $\CC$ is a permutation matrix. Similarly, we get by taking determinants in the constant form of Axiom~\ref{DefModData}.\ref{DefModData4} that
$$\gp^m = (\frac{n_o}{t_o^2})^m \det(\T)^3 \det(\V)$$
Combining both equations and using Proposition~\ref{DefModData}.\ref{DefModData5}, we get that
$$\gp^{2m} = \pm (\frac{n_o}{t_o^2})^{2m} \det(\T)^{6} n^m
= \pm (\frac{\gp\gm}{t_o^4})^{m} \det(\T)^{6}$$
which implies that
$$(\frac{\gp}{\gm})^{m} = \pm \frac{1}{t_o^{4m}} \det(\T)^{6}$$
Now $\det(\T)$ is the product of the diagonal elements~$t_i$,
and therefore an $N$-th root of unity. This implies that
$(\gp/\gm)^{2Nm} = 1$. 

\item 
If the modular datum is integral, $\gp$ and~$\gm$ are contained in~$\Q_N$
by construction. But the only roots of unity contained in~$\Q_N$
are the $N$-th roots of unity and their negatives (cf.~\cite{Wash}, Exerc.~2.3, p.~17). This shows that $(\gp/\gm)^{2N} = 1$, and also that
$(\gp/\gm)^{N} = 1$ if~$N$ is even.
\qed
\end{list}
\end{pf}

\subsection[The central charge]{} \label{CentCharge}
As the name suggests, modular data are closely related to the (homogeneous) modular group~$\SL(2,\Z)$. Recall that the modular group is generated by the two matrices 
$$\gv := \begin{pmatrix} 0 & -1 \\ 1 & 0 \end{pmatrix} 
\qquad \text{and} \qquad
\gt := \begin{pmatrix} 1 & 1 \\ 0 & 1 \end{pmatrix}$$
(cf.~\cite{Apos}, Sec.~2.2, Thm.~2.1, p.~28).
It is easy to see that these matrices satisfy the relations
$$\gv^4 = 1 \qquad (\gt\gv)^3 = \gv^2$$
however, it is a nontrivial result that these are defining relations for the modular group (cf.~\cite{CoxMos}, \S~7.2, p.~85; \cite{FineRosen},  Thm.~3.2.3.2, p.~97; \cite{KasTur}, Thm.~A.2, p.~312; \cite{KleinFricke}, \S~II.9.1, p.~454). The connection to modular data now becomes apparent
if we write Axiom~\ref{DefModData}.\ref{DefModData3} and the 
constant form of Axiom~\ref{DefModData}.\ref{DefModData4} in the form
$$\V^2 = n \CC \qquad \qquad \frac{n_o}{t_o^2} (\T \V)^3 = \gp  \V^2$$
From this, we see that the assignment
$$\SL(2,\Z) \rightarrow \PGL(m,K),~\gv \mapsto \V, \gt \mapsto \T$$
defines a projective representation of the modular group, i.e., a group
homomorphism to the projective general linear group 
$\PGL(m,K):=\GL(m,K)/K^\times$, where~$m$ is the cardinality of~$I$ (cf.~Paragraph~\ref{GroupCohom}).

As we will see in a moment, it follows from the above presentation that every projective representation of the modular group can
be lifted to an ordinary representation, at least if the base field~$K$ is algebraically closed. But doing this in our case requires to
rescale the Verlinde matrix and the Dehn matrix by two scalars, and
as these scalars are not completely determined, certain choices can be
made. We now introduce the notion of an extended modular datum for a modular datum for which these choices have been made:
\begin{defn}
Suppose that $K$ is a field of characteristic zero.
An extended modular datum over~$K$ is a septuple 
$(I,o,{}^*,\V,\T,D,\ell)$ such that $(I,o,{}^*,\V,\T)$
is a modular datum and~$D$ and~$\ell$ are elements of~$K$ that satisfy
$$D^4 = n^2 \qquad \ell^3 = \frac{\gp}{n_o t_o D}$$
where~$n$ is the global dimension and~$\gp$ is the Gaussian sum of the
modular datum. The number~$D$ is called the generalized rank of the extended modular datum. A generalized rank is called a rank if~$D^2=n$.
The number~$\ell$ is called the (multiplicative) central charge of the 
extended modular datum. The matrices
$$\V' := \frac{\V}{D} \qquad \text{and} \qquad \T' := \frac{\T}{t_o \ell}$$
are called the homogeneous Verlinde matrix resp.\ the homogeneous Dehn
matrix of the extended modular datum.
\end{defn}

The notion of rank is adapted from \cite{Tur}, Sec.~II.1.6, p.~76. If we rescale the Verlinde matrix
with a rank, Axiom~\ref{DefModData}.\ref{DefModData3} takes the form
$(\V/D)^2 = \CC$. Note that it follows from Proposition~\ref{DefModData}.\ref{DefModData5}
that, if $\gp=\gm$, the two possible ranks are $D= \pm \gp/n_o$, so that the condition on the central charge becomes $\ell^3= \pm 1/t_o$.

On the other hand, we know from Axiom~\ref{DefModData}.\ref{DefModData3} that the Verlinde matrix has at most four different eigenvalues, namely exactly the four possible values for a generalized rank. If we rescale the Verlinde matrix with a generalized rank, we do not always have the equation $(\V/D)^2 = \CC$, but in any case we have $(\V/D)^4 = \E$, which is one of the defining relations of the modular group. If 
$\gp = \pm \gm$, it follows again from Proposition~\ref{DefModData}.\ref{DefModData5}
that $D=\gp/n_o$ is a possible generalized rank, and the other possibilities deviate from this one by a fourth root of unity.
If $D=\gp/n_o$, the condition on the central charge becomes \mbox{$\ell^3=1/t_o$}.

The following observation is easy but important:
\begin{prop}
The central charge is a root of unity.
\end{prop}
\begin{pf}
From Proposition~\ref{DefModData}.\ref{DefModData5}, we get that
$$\ell^3 = \frac{\gp}{n_o t_o D} = \pm \frac{D n_o}{t_o \gm}$$
which implies that $t_o^2 \ell^6 = \pm \gp/\gm$. Now Proposition~\ref{GaussSum} implies the assertion.~\qed
\end{pf}

Let us now see which role these quantities play in the lifting of the projective representation of the modular group to an ordinary linear representation. With a generalized rank~$D$ and the multiplicative central charge~$\ell$, we can write Axiom~\ref{DefModData}.\ref{DefModData3} and 
Axiom~\ref{DefModData}.\ref{DefModData4} in the form
$$\V^2 = \pm D^2 \CC \qquad \qquad (\T \V)^3 = D t_o^3 \ell^3 \V^2$$
which implies
$$(\frac{\V}{D})^4 = \E \qquad \qquad 
(\frac{\T}{t_o \ell} \frac{\V}{D})^3 = (\frac{\V}{D})^2$$
so that the assignment
$$\SL(2,\Z) \rightarrow \GL(m,K),~\gv \mapsto \V' = \frac{\V}{D}, 
\gt \mapsto \T' = \frac{\T}{t_o \ell}$$
yields a linear representation of the modular group.

We note that extended modular data can be rescaled in the 
following way: If~$\mu$ is a nonzero number and~$\zeta$
is a root of unity, then $(I,o,{}^*,\mu\V,\zeta\T,\mu D,\ell)$
is again an extended modular datum, because the Gaussian sum corresponding to the rescaled matrices is~$\mu^2 \zeta \gp$. Clearly, the homogeneous Verlinde matrix and the homogeneous Dehn matrix are not changed by this rescaling, so that the resulting representation of the modular group is the same. We will discuss still another way of modifying extended modular data in Paragraph~\ref{PropMod}.

It is instructive to consider the case of an integral modular datum over the field~$K=\C$: In this case,
we have $D= \pm \sqrt{n} \in \R$ for a rank, and 
$D= \pm i \sqrt{n}$ for a generalized rank that is not a rank.
It then follows from Proposition~\ref{ComplConj} that the homogeneous Verlinde matrix~$\V'$ becomes unitary, a property that the Dehn matrix and the homogeneous Dehn matrix always have, as their diagonal entries are roots of unity. In this case, we can also write~$\ell$ in the form
$$\ell = e^{2\pi i c/24}$$
for a number $c \in \C$, which is unique modulo~$24\Z$. This number~$c$ is then called the additive central charge of the extended modular datum. As we saw above, $\ell$ is a root of unity, which means that~$c$ is in fact a rational number. 

Our exposition in this paragraph follows essentially~\cite{BakKir}, Sec.~3.1, in particular Eq.~(3.1.17), p.~51 and Rem.~3.1.20, p.~58 (cf.~also \cite{EtNikOst}, Prop.~3.4, p.~595). We note that in \cite{SoZhu} and \cite{Tur}, the conventions are slightly different; essentially, the matrices~$\V$ and~$\V\CC$ are interchanged 
(cf.~\cite{SoZhu}, Prop.~5.3, p.~48; \cite{Tur}, Sec.~II.3.9, p.~98). 
A more comprehensive comparison of the notation used in some of these references can be found in the table at the end of Paragraph~\ref{FactHopf}.

\newpage
\section{Odd Exponents} \label{Sec:OddExp}
\subsection[Properties of modular data]{} \label{PropMod}
We have already introduced the notions of normalization and integrality
for modular data, but we will need some further properties. Recall that the principal congruence subgroup of level~$N$ is the subgroup
$$\Gamma(N) := \{\begin{pmatrix}
a & b \\
c & d
\end{pmatrix} \in \SL(2,\Z) \mid 
a \equiv d \equiv 1, \; b \equiv c \equiv 0 \pmod{N}\}$$
of~$\SL(2,\Z)$; in other words, $\Gamma(N)$ is the kernel of the natural map 
$$\SL(2,\Z) \rightarrow \SL(2,\Z_N)$$
that reduces all entries of a matrix modulo~$N$. 

\begin{defb}
\begin{enumerate}
\item 
A modular datum is called a projective congruence datum if~$\Gamma(N_o)$
is contained in the kernel of the projective representation
$$\SL(2,\Z) \rightarrow \PGL(m,K),~\gv \mapsto \V, \gt \mapsto \T$$
where $N_o$ is the normalized exponent of the modular datum.

\item 
An extended modular datum is called a congruence datum if~$\Gamma(N_o)$ is contained in the kernel of the (ordinary) representation
$$\SL(2,\Z) \rightarrow \GL(m,K),~\gv \mapsto \V', 
\gt \mapsto \T'$$
where~$\V'$ is the homogeneous Verlinde matrix and~$\T'$ is the homogeneous
Dehn matrix.

\item \label{GaloisDat} 
An integral modular datum is called Galois if we have
$$t_{\sigma.i}=\sigma^2(t_i)$$
for all $i \in I$ and all $\sigma \in \Gal(\Q_{N}/\Q)$.
\end{enumerate}
\end{defb}

The notions that we have introduced for modular data of course also apply to extended modular data, and, in this sense, a congruence datum is clearly a projective congruence datum. Furthermore, it is obvious that the property of being a projective congruence datum remains unchanged if we rescale the Verlinde and the Dehn matrix by appropriate factors, because this rescaling does not affect~$N_o$. This also holds for the property of being a congruence datum if we rescale as indicated in Paragraph~\ref{CentCharge}, as we noted there that the homogeneous Verlinde matrix and the homogeneous Dehn matrix remain unchanged.
On the other hand, the property of being Galois is not invariant under
rescaling. We also note that, by definition, every Galois modular datum
is integral.

Even if it is possible to extend a projective congruence datum to 
a congruence datum by choosing a generalized rank~$D$ and a central charge~$\ell$, a different choice of these two numbers will in general not extend this projective congruence datum to a congruence datum. To understand this, let us first consider in how many ways a given modular datum can be extended:
\begin{lemma}
Suppose that $(I,o,{}^*,\V,\T,D,\ell)$ is an extended modular datum.
\begin{enumerate}
\item 
If~$\zeta$ is a twelfth root of unity, then 
$(I,o,{}^*,\V,\T,D/\zeta^3,\zeta\ell)$ is again an extended modular datum.

\item
Conversely, if $(I,o,{}^*,\V,\T,\tilde{D},\tilde{\ell})$ is an extended modular datum, then there exists a twelfth root of unity~$\zeta$
such that~$\tilde{D} = D/\zeta^3$ and $\tilde{\ell}=\zeta\ell$.
\end{enumerate}
\end{lemma}
\begin{pf}
The first assertion follows immediately from Definition~\ref{CentCharge}. For the second assertion, we get from the conditions
$D^4 = n^2 = \tilde{D}^4$ that $D/\tilde{D}$ is a fourth root of unity, and the equations
$$\ell^3 = \frac{\gp}{n_o t_o D} \qquad  \qquad 
\tilde{\ell}^3 = \frac{\gp}{n_o t_o \tilde{D}} $$
imply that~$\tilde{\ell}^3/\ell^3=D/\tilde{D}$, so that~$\tilde{\ell}/\ell$ is a twelfth root of unity. 
\qed
\end{pf}

This lemma tells us that extensions are unique up to twelfth roots of unity. However, if~$D$ was a rank and we also want that~$\tilde{D}$
is a rank, then we are limited to sixth roots of unity. And if we even
want that~$D=\tilde{D}$, then we can only use third roots of unity
for the modification.

Returning to the original question, consider a congruence datum
with generalized rank~$D$ and central charge~$\ell$, and introduce, as in the lemma, a new generalized rank~$\tilde{D} = D/\zeta^3$ and a new central charge~$\tilde{\ell}=\zeta\ell$, where~$\zeta$ is a twelfth root of unity.
If the new choices also lead to a congruence datum, then the two representations 
$$\gv \mapsto \frac{\V}{D}, \; \gt \mapsto \frac{\T}{t_o \ell} 
\qquad \text{and} \qquad 
\gv \mapsto \frac{\V}{\tilde{D}}, \; 
\gt \mapsto \frac{\T}{t_o \tilde{\ell}} $$
of the modular group both factor over the reduced modular group~$\SL(2,\Z_{N_o})$, which implies that the map
$$\SL(2,\Z_{N_o}) \rightarrow K^\times,~\bar{\gv} \mapsto \frac{\tilde{D}}{D}=\frac{1}{\zeta^3},~\bar{\gt} \mapsto \frac{\tilde{\ell}}{\ell}=\zeta$$
defines a one-dimensional representation and therefore
factors over the commutator factor group $\SL(2,\Z_{N_o})_{\ab}$.
But this commutator factor group is cyclic, and its order is the 
greatest common divisor of~${N_o}$ and~$12$ (cf.~\cite{Beyl}, Lem.~(1.13), p.~25). In particular, if~${N_o}$ is an odd number that is not divisible by~$3$, then the commutator factor group is trivial, and we must have~$\zeta=1$, $\tilde{D}=D$, and $\tilde{\ell}=\ell$, so that our generalized rank and our central charge are uniquely determined by the congruence requirement. In Paragraph~\ref{ExamRadf}, we will see explicit examples where~$N_o$ is odd and not divisible by~$3$.

The following table summarizes in short form the properties of modular data that we have introduced:
\begin{center}
\setlength{\extrarowheight}{1pt}
\begin{tabular}[t]{|c|c|c|}\hline 
Property & Definition & Paragraph \\ \hline 
normalized & $n_o = t_o = 1$  & \ref{DefModData} \\ \hline 
integral & $n_i \in \N$  & \ref{DefModData} \\ \hline 
extended & $D$ and $\ell$ chosen  & \ref{CentCharge} \\ \hline 
projective congruence & $\SL(2,\Z_{N_o}) \rightarrow \PGL(m,K)$  & \ref{PropMod} \\ \hline 
congruence & $\SL(2,\Z_{N_o}) \rightarrow \GL(m,K)$  & \ref{PropMod} \\ \hline 
Galois & $t_{\sigma.i}=\sigma^2(t_i)$  & \ref{PropMod} \\ \hline 
\end{tabular} 
\end{center}

\subsection[The fusion symbol]{} \label{FusSym}
Recall the classical isomorphism
$$\Z_N^\times \rightarrow \Gal(\Q_N/\Q),~\bar{q} \mapsto \sigma_q$$
between the group of units $\Z_N^\times$ of the ring~$\Z_N$
and the Galois group of the cyclotomic field that maps a residue class~$\bar{q}$ to the automorphism~$\sigma_q$ that raises every $N$-th root of unity to its $q$-th power. By analogy with the Hopf symbol introduced in \cite{SoZhu}, Par.~12.1, p.~114, we associate with every integral modular datum its fusion symbol:
\begin{defn}
Consider an integral modular datum with exponent~$N$ and Gaussian sum~$\gp$.
For $q \in \Z$, we define the fusion symbol
$$
\fs(q) := \begin{cases} \displaystyle\frac{\sigma_q(\gp)}{\gp}  &: \gcd(q,N) = 1 \\
0 &: \gcd(q,N) \neq 1
\end{cases} $$
\end{defn}
The fusion symbol should be viewed as a 1-cocycle
in the following way (cf.~\cite{Serre2}, Chap.~VII, \S~3, p.~113):
The Galois group~$\Gal(\Q_N/\Q)$ acts on
the multiplicative group~$\Q_N^\times$ of nonzero elements in the cyclotomic field. The element $\gp \in \Q_N^\times$ then can be viewed as a 0-cochain, and therefore gives rise to a 1-coboundary
$$\Gal(\Q_N/\Q) \rightarrow \Q_N^\times,~\sigma 
\mapsto \frac{\sigma(\gp)}{\gp}$$
which is essentially the fusion symbol, as~$\sigma_q$ is mapped to~$\fs(q)$. Now every \mbox{1-coboundary} is in particular a 1-cocycle, and the 1-cocycle condition means for the fusion symbol that
$$\fs(qq') = \fs(q) \sigma_q(\fs(q'))$$
Actually, it follows from a version of Hilbert's theorem~90 that every cocycle is a coboundary in this situation (cf.~\cite{Serre2}, Chap.~X, \S~1, Prop.~2, p.~150). This equation also shows that the fusion symbol is a Dirichlet character, i.e., satisfies
$$\fs(qq') = \fs(q) \fs(q')$$
for all $q,q' \in \Z_N^\times$, if and only if 
$\sigma_q(\fs(q')) = \fs(q')$, i.e., if and only if~$\fs(q')$
is invariant under the Galois group, which means that it is a rational number. In this case, we have $\fs(q')=\pm 1$, because $\{1,-1\}$
is the largest finite multiplicative subgroup of~$\Q^\times$.

It is obvious that~$\fs(q)$ depends only on the residue class of~$q$ modulo~$N$, and also that $\fs(1)=1$ and $\fs(-1)=\gm/\gp$. Therefore, the following result generalizes Proposition~\ref{GaussSum}:
\begin{prop}
Suppose that the modular datum is integral, and that $q \in \Z$ is relatively prime to~$N$. Then we have $\fs(q)^{2N}=1$. If $N$ is even, we have $\fs(q)^{N}=1$.
\end{prop}
\begin{pf}
We have already seen in Paragraph~\ref{GaussSum} that
$\det(\V)^2 = \pm n^m$, where~$m$ is the cardinality of~$I$. Since $\det(\V) \in \Q_N$ by Proposition~\ref{DefModData}.\ref{DefModDataST}, this equation implies
that $\sigma_q(\det(\V)) = \pm \det(\V)$. We have also seen in Paragraph~\ref{GaussSum} that
$$\gp^m = (\frac{n_o}{t_o^2})^m \det(\T)^3 \det(\V)$$
Applying~$\sigma_q$ to this equation, we get
$$\sigma_q(\gp)^m = 
\pm (\frac{n_o}{t_o^{2q}})^m \sigma_q(\det(\T))^3 \det(\V)$$
If we divide the two equations by each other, we therefore get that
$$\fs(q)^m = \frac{\sigma_q(\gp)^m}{\gp^m} = 
\pm (\frac{t_{o}^2}{t_o^{2q}})^m \frac{\sigma_q(\det(\T))^3}{\det(\T)^3}$$
Since $\det(\T)$ and~$t_o$ are $N$-th roots of unity, this implies
that $\fs(q)^{2Nm} = 1$. But $\fs(q) \in \Q_N$ by definition, and if
$N$ is even, the only roots of unity contained in this cyclotomic field
are the $N$-th roots of unity (cf.~\cite{Wash}, Exerc.~2.3, p.~17). If $N$
is odd, it still follows from the same argument that at least
$\pm \fs(q)$ is an $N$-th root of unity, which implies the assertion.
\qed
\end{pf}

One of our main goals in the sequel is to sharpen the preceding proposition. As a first step, we record the following consequence, continuing to assume that the modular datum is integral:
\begin{corollary}
The fusion symbol is a Dirichlet character if and only if $\gm = \pm \gp$.
In this case, we have $\fs(q) \in \{0,1,-1\}$ for all $q \in \Z$.
\end{corollary}
\begin{pf}
We have already discussed above that $\fs(q) \in \{0,1,-1\}$
if the fusion symbol is a Dirichlet character, and also that in this
case $\gm/\gp = \fs(-1) = \pm 1$. So, suppose now that conversely
$\gm = \pm \gp$. Expressed differently, this means that 
$\gamma(\sum_{i \in I} n_i^2 t_i) = \pm \sum_{i \in I} n_i^2 t_i$.
If~$q$ is relatively prime to~$N$, we can apply the automorphism~$\sigma_q$
to this equation to get
$\gamma(\sum_{i \in I} n_i^2 t_i^q) = \pm \sum_{i \in I} n_i^2 t_i^q$,
where we have used that~$\sigma_q$ commutes with~$\gamma$ because the Galois group is abelian. Dividing this equation by the preceding one, we obtain $\gamma(\fs(q)) = \fs(q)$. But since~$\fs(q)$ is a root of unity by the preceding result, we also have $\gamma(\fs(q)) = 1/\fs(q)$. Therefore $\fs(q)^2=1$ and $\fs(q) \in \{1,-1\}$. 
\qed
\end{pf}

\subsection[Twelfth roots of unity]{} \label{12root}
Our second step to improve Proposition~\ref{FusSym} is to show that, for Galois modular data, the fusion symbols are in fact twelfth roots of unity. For this, we recall the so-called `definition of~24' (cf. \cite{CostGann}, \S~2.4, Prop.~3.b, p.~9; \cite{CostGannRue}, Sec.~2.2, p.~691; \cite{GanMoonMonst}, Sec.~2.5.1, p.~168):
\begin{lemma}
Suppose that~$\xi$ is an~$N$-th root of unity, and that $\xi^{q^2}=\xi$ for all~$q \in \Z$ that 
are relatively prime to~$N$. Then $\xi^{24}=1$. 
\end{lemma}
\begin{pf}
Recall that our base field~$K$ has characteristic zero. Choose a primitive $N$-th root of unity~$\zeta$ in its algebraic closure and write 
$\xi = \zeta^k$ for some~$k$, so that we have
$\zeta^{k q^2} = \zeta^k$. This implies 
$k q^2 \equiv k \pmod{N}$ for all~$q$ that are relatively prime to~$N$, which means that~$N$ divides $k(q^2-1)$. If $N$ is odd, we see by taking~$q=2$ that~$N$ divides~$3k$, so that
$\xi^3 = \zeta^{3k} = 1$.

If~$N$ is even, and $N=\prod_i p_i^{k_i}$ is the prime factorization of~$N$ into powers of distinct primes, we can by the Chinese remainder theorem find a unit~$q$ modulo~$N$ that satisfies $q \equiv 3 \pmod{p_i^{k_i}}$ if~$p_i=2$ and  
$q \equiv 2 \pmod{p_i^{k_i}}$ if~$p_i \neq 2$. If~$p_i=2$, we therefore get that~$p_i^{k_i}$ divides~$8k$, and if~$p_i \neq 2$, $p_i^{k_i}$ divides~$3k$. In any case, $p_i^{k_i}$ divides~$24k$, so that~$N$
divides~$24k$, showing that $\xi^{24} = \zeta^{24k} = 1$.~\qed
\end{pf}

For example, this lemma implies that for Galois modular data we have
\mbox{$t_o^{24}=1$}, as the Galois condition and Proposition~\ref{GaloisAct} imply
$\sigma^2(t_o) = t_{\sigma.o} = t_o$ for all 
\mbox{$\sigma \in \Gal(\Q_N/\Q)$}. It can be used in a similar way to derive our improved statement:
\begin{prop}
For a Galois modular datum, we have $\fs(q)^{12}=1$ for all $q \in \Z$ that are relatively prime to~$N$. 
\end{prop}
\begin{pf}
\begin{list}{(\arabic{num})}{\usecounter{num} \leftmargin0cm \itemindent5pt}
\item
By possibly enlarging the base field, we can choose a rank~$D$ for the modular datum. From Proposition~\ref{DefModData}.\ref{DefModDataRecGauss},
we then have $\gp/n_o D = n_o D/\gm$ and therefore
$$(\frac{\gp}{n_o D})^2 = \frac{\gp}{\gm}$$
so that Proposition~\ref{GaussSum} implies that $\gp/n_oD$
is a $4N$-th root of unity. Because $\gp \in \Q_{N}$ and 
$n_o \in \Z$, this also implies that $D \in \Q_{4N}$.

\item
For any $\sigma \in \Gal(\Q_N/\Q)$, it follows directly from the Galois property and Proposition~\ref{GaloisAct} that 
\begin{align*}
\sigma^2(\gp) = \sum_{i \in I} n_i^2 \sigma^2(t_i) = 
\sum_{i \in I} n_{\sigma.i}^2 t_{\sigma.i} = \gp
\end{align*}
Similarly, we have for $\sigma \in \Gal(\Q_{4N}/\Q)$ that
$\sigma(D)^2 = \sigma(n) = n$ and therefore $\sigma(D) = \pm D$,
which implies that $\sigma^2(D) = D$. Both facts together show that
$\sigma^2(\gp/n_o D) = \gp/n_o D$ for all $\sigma \in \Gal(\Q_{4N}/\Q)$.
Now the lemma above implies that $(\gp/n_o D)^{24}=1$.

\item
Now suppose that $N$ is odd. Since we have already seen that 
$(\gp/n_o D)^{4N}=1$, we then have $(\gp/n_o D)^{12}=1$,
as the greatest common divisor of~$24$ and~$4N$ now divides~$12$.
If~$q$ is relatively prime to~$4N$, we can apply~$\sigma_q$ to this equation to get $(\sigma_q(\gp/n_o D))^{12}=1$, and if we divide
this equation by the preceding one, we get
$$\fs(q)^{12} (\frac{D}{\sigma_q(D)})^{12} =
(\frac{\sigma_q(\gp)}{\gp})^{12} (\frac{D}{\sigma_q(D)})^{12} = 1$$
Now we saw in the preceding step that 
$\sigma_q(D)^2 = n = D^2$, so that our assertion $\fs(q)^{12} = 1$
holds if~$q$ is relatively prime to~$4N$. But as~$\fs(q)$
depends only on the residue class of~$q$ modulo~$N$ and the
canonical map from~$\Z_{4N}^\times$ to~$\Z_{N}^\times$ is surjective, the assertion also holds if~$q$ is only relatively prime to~$N$.

\item
Now suppose that $N$ is even. If~$q$ is relatively prime to~$N$, then
$q$ is odd and also relatively prime to~$4N$. Since $\gp/n_o D$
is a $4N$-th root of unity, we have 
$$\sigma_q(\frac{\gp}{n_o D}) = (\frac{\gp}{n_o D})^q$$
which together with $\sigma_q(D) = \pm D$ implies that
$$\fs(q)= \frac{\sigma_q(\gp)}{\gp} = \pm (\frac{\gp}{n_o D})^{q-1}$$
Because $q$ is odd, $q-1$ is even and $12(q-1)$ is divisible by~$24$, so that we get $\fs(q)^{12} = (\gp/n_o D)^{12(q-1)} = 1$, as asserted.
\qed
\end{list}
\end{pf}

Note that this result means in particular that $\gm/\gp = \fs(-1)$ is a 12-th root of unity. In the case where~$N$ is odd, the preceding result
will be improved in Theorem~\ref{OddExp}.

The lemma above has another consequence that is worth noting:
As we have seen, it follows from Proposition~\ref{DefModData}.\ref{DefModData4} that, for integral
modular data, the entries of the Verlinde matrix are contained in the
cyclotomic field~$\Q_N$, which implies that we have a field
extension 
$$L:=\Q(\{s_{ij}| i,j \in I\}) \subset \Q_N = \Q(\{t_i | i \in I\})$$
about which we can say the following:
\begin{corollary}
If the modular datum is Galois, then the index $[\Q_N:L]$ is a power of~$2$. If in addition~$L=\Q$, then $[\Q_N:L]$ divides~$8$ and~$N$ divides~$24$.
\end{corollary}
\begin{pf}
By the fundamental theorem of Galois theory, the index $[\Q_N:L]$ is equal to the order of the Galois group $\Gal(\Q_N/L)$, which is a subgroup of $\Gal(\Q_N/\Q)$. If $\sigma \in \Gal(\Q_N/L)$, then we have $s_{ij} = \sigma(s_{ij}) = s_{\sigma.i,j}$ by Proposition~\ref{GaloisAct}, which implies that $\sigma.i = i$ for all $i \in I$. By the Galois condition,
we then have $\sigma^2(t_i) = t_{\sigma.i} = t_i$, so that
$\sigma^2=\id$. This shows that $\Gal(\Q_N/L)$ has exponent~$2$.
By Cauchy's theorem, this implies that the order of $\Gal(\Q_N/L)$
is a power of~$2$.

If~$L=\Q$, then this argument shows that $\sigma^2(t_i) = t_i$ for all $\sigma \in \Gal(\Q_N/\Q)$, so that~$t_i^{24}=1$ by the preceding lemma,
which implies that~$N$ divides~$24$. This in turn implies that
$[\Q_N:\Q]$ divides $[\Q_{24}:\Q]=8$.
\qed
\end{pf}

We note that at least the second half of this corollary is present 
in \cite{CostGannRue}, Sec.~3.3, p.~698, which also contains
an example that shows that~$L=\Q$ actually can occur.

\subsection[The relation of the actions]{} \label{RelAct}
Like every permutation, the permutation~$i \mapsto \sigma.i$ that an automorphism $\sigma \in \Gal(\Q_{N_o}/\Q)$ induces on the index set of an integral modular datum gives rise to a permutation matrix, which we denote by
$$P(\sigma) := (\delta_{i,\sigma.j})_{i,j \in I}$$
Note that $P(\gamma) = \CC$ by Proposition~\ref{ComplConj}, so that in particular $P(\sigma)$ commutes with~$\CC$. With the Verlinde matrix, this permutation matrix has the following commutation relation (cf.~\cite{SoZhu}, Cor.~10.1, p.~97):
\begin{lemma}
$\V P(\sigma) = P(\sigma)^{-1} \V$
\end{lemma}
\begin{pf}
From Proposition~\ref{GaloisAct}, we get
\begin{align*}
P(\sigma)^{-1} \V = (\sum_{k \in I} \delta_{i,\sigma^{-1}.k} s_{kj})_{i,j}
= (s_{\sigma.i,j})_{i,j} = (\sigma(s_{i,j}))_{i,j}
= (s_{i,\sigma.j})_{i,j} = \V P(\sigma)
\end{align*}
as asserted.
\qed
\end{pf}

The following argument, which is an adaption of the one found in \cite{CostGann}, \S~2.3, Thm.~2, p.~7f, shows that for Galois modular data
the permutation matrix can be computed from the Verlinde matrix and the Dehn matrix:
\begin{prop}
Consider a Galois modular datum with Verlinde matrix~$\V$ and Dehn matrix~$\T$. Suppose that~$q, q' \in \Z$ satisfy $qq' \equiv 1 \pmod{N}$. 
Then we have
$$\V \T^{q'} \V^{-1} \T^q \V \T^{q'} = 
\frac{t_o^{2q}}{n_o} \sigma_q(\gp)  \; P(\sigma_q^{-1})$$
\end{prop}
\begin{pf}
In the constant form of Axiom~\ref{DefModData}.\ref{DefModData4}, we bring all matrices to one side and write out the resulting equation 
$(\V\T)^3 = \frac{t_o^2}{n_o} \gp n \CC$ in components. Then we get
$$\sum_{j,k \in I} s_{ij} t_j s_{jk} t_k s_{kl} t_l
= \frac{t_o^2}{n_o} \gp n \delta_{i,l^*}$$
If we apply~$\sigma_q$ to this equation, it becomes 
$$\sum_{j,k \in I} s_{i,\sigma_q.j} t_j^q s_{\sigma_q.j,k} t_k^q s_{k,\sigma_q.l} t_l^q = \frac{n}{n_o} t_o^{2q} \sigma_q(\gp)  \delta_{i,l^*}$$
by Proposition~\ref{GaloisAct}.
If we replace~$j$ by~$\sigma_q^{-1}.j$ and~$l$ by~$\sigma_q^{-1}.l$,
this becomes 
$$\frac{n}{n_o} t_o^{2q} \sigma_q(\gp) \delta_{i,\sigma_q^{-1}.l^*} = \sum_{j,k \in I} s_{ij} t_{\sigma_q^{-1}.j}^q s_{jk} t_k^q s_{kl} t_{\sigma_q^{-1}.l}^q
= \sum_{j,k \in I} s_{ij} t_{j}^{q'} s_{jk} t_k^q s_{kl} t_{l}^{q'}$$
where in the last step we have used the Galois hypothesis. But this is the component form of the matrix equation
$$\V \T^{q'} \V \T^q \V \T^{q'} = 
\frac{n}{n_o} t_o^{2q} \sigma_q(\gp)  \; P(\sigma_q^{-1}) \CC$$
which is equivalent to the assertion.
\qed
\end{pf}

\subsection[Odd exponents]{} \label{OddExp}
From the relation between the Verlinde matrix, the Dehn matrix, and the permutation matrices given in Proposition~\ref{RelAct}, we can derive
the following relation for the matrix elements:
\begin{prop}
Consider a Galois modular datum, and suppose that~$q$ and~$r$ are relatively prime to~$N$. Choose numbers~$q'$ and~$r'$ such that 
$qq' \equiv rr' \equiv 1 \pmod{N}$.
Then we have 
\begin{align*}
&\sigma_{r}(\gp) \sum_{k \in I} s_{i^*,k} s_{jk} t_k^{q-r} = 
t_o^{2(q-r)} \sigma_{q}(\gp) \; t_i^{-q'} t_j^{r'}  
\sum_{k \in I} s_{\sigma_{q'}.i,k^*} s_{\sigma_{r'}.j,k} t_k^{r'-q'}
\end{align*}
\end{prop}
\begin{pf}
From Proposition~\ref{RelAct}, it follows that 
$$\V^{-1} \T^q \V = 
\frac{t_o^{2q}}{n_o} \sigma_q(\gp) \; 
\T^{-q'} \V^{-1} P(\sigma_q^{-1}) \T^{-q'}$$
Inverting this relation and substituting~$r$ for~$q$, we get
$$\sigma_r(\gp) \V^{-1} \T^{-r} \V = 
\frac{n_o}{t_o^{2r}} \; \T^{r'} P(\sigma_r) \V  \T^{r'}$$
Multiplying these two equations, we get
\begin{align*}
\sigma_r(\gp) \V^{-1} \T^{q-r} \V =
t_o^{2(q-r)} \sigma_q(\gp) \; 
\T^{-q'} \V^{-1} P(\sigma_q^{-1}) \T^{r'-q'} P(\sigma_r) \V  \T^{r'}
\end{align*}
If we multiply this equation by the global dimension~$n$ and use
that $n \V^{-1} = \V\CC$, we arrive at the assertion in matrix form.
\qed
\end{pf}

The case $q=-1$, $r=1$ of this proposition holds without the Galois condition:
\begin{corollary}
For all $i,j \in I$, we have
\begin{align*}
&\gp \sum_{k \in I} s_{i^*,k} s_{jk} t_k^{-2} = 
\gm \frac{t_i t_j}{t_o^4} \sum_{k\in I} s_{ik} s_{jk} t_k^{2}
\end{align*}
\end{corollary}
\begin{pf}
By squaring the constant form of Axiom~\ref{DefModData}.\ref{DefModData4},
we get
$$\gp^2 (\T^{-1} \V \T^{-1})^2 = \frac{n_o^2}{t_o^4} (\V \T \V)^2 = \frac{\gp \gm}{n t_o^4} (\V \T \V)^2$$
where we have used Proposition~\ref{DefModData}.\ref{DefModDataRecGauss}
for the second equality. Cancelling one Gaussian sum and using Axiom~\ref{DefModData}.\ref{DefModData3}, this becomes
$$\gp \T^{-1} \V \T^{-2} \V \T^{-1} = \frac{\gm}{t_o^4} \V \T \CC \T \V$$
Bringing the outer Dehn matrices to the right-hand side and the
charge conjugation matrix to the left-hand side, this equation can be written in the equivalent form
$$\gp \CC \V \T^{-2} \V = \frac{\gm}{t_o^4} \T \V \T^2 \V \T$$
The assertion is the $(i,j)$-component of this matrix equation.
\qed
\end{pf}

From this, we can now take the third step to sharpen Proposition~\ref{FusSym}. We have considered the case of odd exponents already in the proof of Proposition~\ref{12root}. As we will see now, in this case the fusion symbol is not only a twelfth, but rather a second root of unity:
\begin{thm}
If the exponent~$N$ of an integral modular datum is odd, we have
$\gp = \pm t_o^2 \gm$. If the modular datum is in addition normalized, we have $\gp = \pm \gm$ and $\fs(q)= \pm 1$ if $q$ is relatively prime to~$N$.
\end{thm}
\begin{pf}
\begin{list}{(\arabic{num})}{\usecounter{num} \leftmargin0cm \itemindent5pt}
\item
By Proposition~\ref{GaussSum}, $\gm/\gp$
is a $2N$-th root of unity. We can therefore write it in the form
$\gm/\gp = v \zeta$, where $v \in \{1,-1\}$ and~$\zeta$ is an $N$-th root of unity, so that $\gm = v \zeta \gp$.
Because~$N$ is odd, we can apply the Galois automorphism $\sigma_2 \in \Gal(\Q_N/\Q)$ to this formula. If we also divide by~$\gp$, we get
$$\fs(-2) = \frac{\sigma_{-2}(\gp)}{\gp} = 
\frac{1}{\gp} \sigma_2(\gm) = \frac{1}{\gp} v \zeta^2 \sigma_2(\gp) = 
v \zeta^2 \fs(2)$$

\item
On the other hand, setting $i=j=o$ in the preceding corollary, we get
\begin{align*}
&\gp \sum_{k \in I} n_{k}^2 t_k^{-2} = 
\frac{\gm}{t_o^2} \sum_{k\in I} n_{k}^2 t_k^{2}
\end{align*}
If we divide by~$\gp^2$, this formula becomes 
$$\fs(-2) = \frac{v \zeta}{t_o^2} \fs(2)$$ Comparing this with the formula obtained in the first step, we see that~$\zeta=1/t_o^2$, which establishes the first assertion.

\item
For a normalized modular datum, we have $t_o=1$ by definition, so that
$\gp= \pm \gm$. By Corollary~\ref{FusSym}, this implies $\fs(q) = \pm 1$
if~$q$ is relatively prime to~$N$.~\qed
\end{list}
\end{pf}

For an extended integral modular datum with odd exponent, the preceding theorem means for the multiplicative central charge that~$\ell^{12}=1$, because we have seen in the proof of Proposition~\ref{CentCharge} that $\ell^6 = \pm \gp/(t_o^2 \gm)$.
In the case $K=\C$, we have also explained in Paragraph~\ref{CentCharge} that we can write $\ell = e^{2\pi i c/24}$, where~$c$ is the additive central charge. The equation $\ell^{12}=1$ then means that~$c$ must be an even integer. In Paragraph~\ref{Semion}, we will see in an explicit example that the assumption that the exponent is odd cannot be omitted.

\subsection[The classical Gaussian sum]{} \label{ClassGauss}
So far, it has not been necessary to consider classical Gaussian sums,
and it is not yet apparent from our considerations why the quantity~$\gp$ is called the Gaussian sum of the modular datum. We will explain this in Paragraph~\ref{ExamRadf}; however, as we will see now, a comparison between the Gaussian sum of a modular datum and the classical Gaussian sum leads to important new insights. Recall the 
definition of the classical Gaussian sum (cf.~\cite{Nag}, Chap.~V, Sec.~53, p.~177): 
\begin{defn}
For a primitive $n$-th root of unity~$\zeta \in K$, the classical 
Gaussian sum~$\Gp$ is defined as 
$$\Gp = \Gp_n(\zeta) := \sum_{i=0}^{n-1} \zeta^{i^2}$$
We also define $\Gm = \Gp_n(1/\zeta)$.
\end{defn}

The basic facts about classical Gaussian sums are summarized in the
following nontrivial lemma, which is mainly due to C.~F.~Gau\ss{} (cf.~\cite{GaussSum}):
\begin{lemma}
If $\iota \in K$ is a primitive fourth root of unity, we have
$$\Gp^2 = \begin{cases}
\pm 2 \iota n &: n \equiv 0 \pmod{4} \\
n &: n \equiv 1 \pmod{4} \\
0 &: n \equiv 2 \pmod{4} \\
-n &: n \equiv 3 \pmod{4} 
  \end{cases}$$
If $n$ is odd and~$q$ is relatively prime to~$n$, we have
$$\sigma_q(\Gp) = \Gp_n(\zeta^q) = \jac{q}{n} \Gp$$
where $\jac{q}{n}$ is the Jacobi symbol.
\end{lemma}
\begin{pf}
The first statement is proved in \cite{Nag}, Chap.~V, Sec.~53, Thm.~99, p.~177, see also~\cite{BeEvWil}, Sec.~1.2, Cor.~1.2.3, p.~15 and~\cite{Lan}, Chap.~IV.VI, Thm.~211, p.~197. The second statement 
is proved in \cite{RadElemNum}, Chap.~11, Thm.~38, p.~93, see also~\cite{BeEvWil}, Sec.~1.5, Thm.~1.5.2, p.~26 and 
\cite{Nag}, Chap.~V, Exerc.~122, p.~187.
\qed
\end{pf}

Under any embedding of the cyclotomic field~$\Q_n \subset K$
into the complex numbers~$\C$, $\Gm$ becomes the complex conjugate
of~$\Gp$, and therefore~$\Gp \Gm$ becomes the square of the absolute value
of~$\Gp$. Using this, we get from the lemma that
$$\Gp \Gm = \begin{cases}
2  n &: n \equiv 0 \pmod{4} \\
0 &: n \equiv 2 \pmod{4}  \\
n &: n \equiv 1 \; \text{or} \; 3 \pmod{4}
\end{cases}$$

In the case where also~$n$ is odd, we can use these classical results to
relate the fusion symbol to the Jacobi symbol:
\begin{thm}
Suppose that both the exponent~$N$ and the global dimension~$n$ of an integral modular datum are odd.
Then we have 
$$t_o^2 \gm = (-1)^{\frac{n-1}{2}} \gp = 
\begin{cases}
\gp &: n \equiv 1 \pmod{4} \\
- \gp &: n \equiv  3 \pmod{4}
\end{cases}$$
If the modular datum is in addition normalized, we have 
$\fs(q) = \jac{q}{n}$ for all $q \in \Z$ that are relatively prime to~$Nn$.
\end{thm}
\begin{pf}
We can assume that~$K$ is algebraically closed, so that we can
choose a primitive $n$-th root of unity~$\zeta$ to form the classical 
Gaussian sum $\Gp=\Gp_n(\zeta)$. From Theorem~\ref{OddExp}, we know that there is a sign~$v \in \{1,-1\}$
such that $\gp = v t_o^2 \gm$. Furthermore, we have $\gp \gm = n n_o^2$ by Proposition~\ref{DefModData}.\ref{DefModDataRecGauss}, and therefore
$$(\frac{\gp}{t_o n_o})^2 = \frac{v \gp \gm}{n_o^2} = v n$$
If $v \neq (-1)^{\frac{n-1}{2}}$, i.e., $v = - (-1)^{\frac{n-1}{2}}$,
we have by the preceding lemma that 
$\Gp^2 = (-1)^{\frac{n-1}{2}} n = -v n$
and therefore
$$(\frac{\gp}{t_o n_o \Gp})^2 = -1$$
showing that $\gp/(t_o n_o \Gp)$ is a primitive fourth root of unity.
But as we have $\gp/(t_o n_o) \in \Q_N$ and $\Gp \in \Q_n$ by construction, we also have $\gp/(t_o n_o \Gp) \in \Q_{Nn}$. As $Nn$ is odd, this implies by~\cite{Wash}, Exerc.~2.3, p.~17 that
$$1= (\frac{\gp}{t_o n_o \Gp})^{2Nn} = (-1)^{Nn} = -1$$
a contradiction. We therefore must have
$v = (-1)^{\frac{n-1}{2}}$, which on the one hand proves the first assertion, on the other hand shows that 
$(\frac{\gp}{t_o n_o})^2 = \Gp^2$,
in other words, that $\gp/(t_o n_o) = \pm \Gp$.

If the modular datum is in addition normalized, this means that
$\gm = (-1)^{\frac{n-1}{2}} \gp$ and $\gp = \pm \Gp$. If~$q$ is relatively prime to~$Nn$, this implies that
$$\fs(q) \gp = \sigma_q(\gp) = \pm \sigma_q(\Gp) = 
\pm \jac{q}{n} \Gp = \jac{q}{n} \gp$$
where the first equation holds by the definition of the fusion symbol and the third by the preceding lemma. This proves the second assertion.
\qed
\end{pf}

The conclusion of this theorem cannot be reformulated in terms of the central charge without the assumption that the generalized rank is actually a rank. So, suppose that we extend the modular datum considered in the preceding theorem by choosing a rank~$D$ and a multiplicative central charge~$\ell$. We then get directly from Definition~\ref{CentCharge} and Proposition~\ref{DefModData}.\ref{DefModDataRecGauss} that
$$\ell^6 = (\frac{\gp}{n_o t_o D})^2 
= (-1)^{\frac{n-1}{2}} \frac{\gp \gm}{n_o^2 n} = (-1)^{\frac{n-1}{2}}$$
In the case $K=\C$, we write, as explained in Paragraph~\ref{CentCharge},
$\ell = e^{2\pi i c/24}$, and then find the condition
\begin{align*}
&c \equiv 0 \pmod{4} \quad \text{if} \quad n \equiv 1 \pmod{4} \\
&c \equiv 2 \pmod{4} \quad \text{if} \quad n \equiv 3 \pmod{4} 
\end{align*}
for the additive central charge~$c$.

The hypothesis of this theorem raises the question how the prime factors of~$n$ and~$N$ are related. We will consider this question again in Paragraph~\ref{Cauchy}.

\newpage
\section{Even Exponents} \label{Sec:EvenExp}
\subsection[Group cohomology]{} \label{GroupCohom}
We will need some facts from group cohomology. Suppose that~$G$ is a finite group and that~$K$ is an algebraically closed field of characteristic zero. Recall that
we have already encountered 1-cocycles in Paragraph~\ref{FusSym}. 
A 2-cocycle of~$G$ with values in the multiplicative group~$K^\times$,
endowed with the trivial $G$-action, is a function
$\omega: G \times G \rightarrow K^\times$ that satisfies
$$\omega(g,g' g'') \omega(g',g'') = 
\omega(g g',g'') \omega(g,g')$$
for all $g, g', g'' \in G$. The set~$Z^2(G,K^\times)$ of all
2-cocycles is a group under pointwise multiplication.
For any map $\nu: G \rightarrow K^\times$, the equation
$$\omega(g,g') := \frac{\nu(g) \nu(g')}{\nu(g g')}$$
defines a 2-cocycle; 2-cocycles of this form are called
2-coboundaries and form a subgroup that is denoted by~$B^2(G,K^\times)$.
The quotient group $H^2(G,K^\times) := Z^2(G,K^\times)/B^2(G,K^\times)$
is called the second cohomology group of~$G$ with coefficients in~$K^\times$, or the Schur multiplicator of~$G$. This group is obviously abelian; it can be shown that it is also finite (cf.~\cite{Hup1}, Kap.~V, Hilfssatz~23.2, p.~629). Moreover, the Schur multiplicators arising from different fields are isomorphic, as long as these fields satisfy our
hypotheses above (cf.~\cite{Hup1}, Kap.~V, Hauptsatz~23.5, p.~631).

The Schur multiplicator is closely related to projective representations,
as explained in \cite{Hup1}, Kap.~V, \S~24, p.~638ff. We will need a slight
refinement of the discussion there. As usual, we embed the multiplicative
group~$K^\times$ of the base field into~$\GL(m,K)$ by mapping 
every nonzero number to the corresponding multiple of the identity
matrix. The corresponding factor group $\PGL(m,K):=\GL(m,K)/K^\times$
is called the projective linear group. Similarly, if $K^\times_l$ denotes
the subgroup of~$K^\times$ consisting of the $l$-th roots of unity,
we denote the corresponding factor group by $\Pn_l\GL(m,K):=\GL(m,K)/K_l^\times$, so that we have
a canonical group homomorphism from~$\Pn_l\GL(m,K)$ to~$\PGL(m,K)$.
The exponent of the Schur multiplicator now determines how far we can lift
a projective representation:
\begin{prop}
Suppose that $\alpha: G \rightarrow \PGL(m,K)$ is a group
homomorphism, and let~$l$ be the exponent of the Schur multiplicator of~$G$. Then there exists a group homomorphism
$\beta: G \rightarrow \Pn_l \GL(m,K)$ such that the diagram
$$\Atriangle[G`\Pn_l\GL(m,K)`\PGL(m,K);\beta`\alpha`]$$
is commutative.
\end{prop}
\begin{pf}
This follows from a minor modification of a standard argument in
the theory of projective representations (cf.~\cite{Hup1}, Kap.~V, Hilfssatz~23.2, p.~629): 
For an invertible matrix $A \in \GL(m,K)$, we denote its coset in~$\PGL(m,K)$ by~$A K^\times$ and its coset in~$\Pn_l\GL(m,K)$
by~$A K_l^\times$. For every $g \in G$, we choose 
a representative $A(g) \in \GL(m,K)$ such that $\alpha(g) = A(g) K^\times$.
Because~$\alpha$ is a group homomorphism, the matrix
$A(g) A(g') A(gg')^{-1}$ is proportional to the unit matrix.
The corresponding proportionality factor~$\omega(g,g')$ is a 2-cocycle (cf.~\cite{Hup1}, Kap.~V, Hilfssatz~24.2, p.~638). By assumption,
$\omega(g,g')^l$ is a coboundary, so there exists a 1-cochain~$\nu$ such that
$$\omega(g,g')^l =  \frac{\nu(g) \nu(g')}{\nu(gg')}$$
Because~$K$ is algebraically closed, we can choose 
for every $g \in G$ a number $\kappa(g) \in K^\times$ such that $\nu(g) = \kappa(g)^l$.
Then $\omega(g,g')$ and $\kappa(g) \kappa(g')/\kappa(gg')$ differ only by an $l$-th root of unity, which implies that
$$A(g) A(g') \qquad \text{and} \qquad 
\frac{\kappa(g) \kappa(g')}{\kappa(gg')} A(gg')$$
are contained in the same coset of~$K_l^\times \subset \GL(m,K)$.
In other words, if we define $B(g):= \frac{1}{\kappa(g)} A(g)$,
then the map $\beta(g) := B(g) K_l^\times$
satisfies our requirements.
\qed
\end{pf}

We are here interested in the case where the finite group~$G$ is the
reduced modular group~$\SL(2,\Z_N)$. Its Schur multiplicator has
been determined by F.~R.~Beyl (cf.~\cite{Beyl}, Thm.~(3.9), p.~32):
\begin{thm}[F.~R.~Beyl]
$$H^2(\SL(2,\Z_N),K^\times) \cong 
\begin{cases}
\Z_2 &: 4 \mid N \\
\{0\} &: 4 \nmid N
\end{cases}$$
\end{thm}

If we combine this theorem with the preceding proposition, we get
the following corollary, in which we have formulated the conclusion
in a slightly different way:
\begin{corollary}
Suppose that $\alpha: \SL(2,\Z_N) \rightarrow \PGL(m,K)$ is a group
homomorphism. Then we can choose, for every $g \in \SL(2,\Z_N)$, 
a representative~$A(g)$ of the coset~$\alpha(g)$ in such a way that
$$A(g) A(g') = \pm A(gg')$$
for all $g, g' \in \SL(2,\Z_N)$.
\end{corollary}
\begin{pf}
As the theorem shows, the exponent~$l$ of the Schur multiplicator
divides~$2$. As we have $K_2^\times = \{1, -1\}$, this implies the assertion.
\qed
\end{pf}

\subsection[Diagonal matrices]{} \label{DiagMat}
To proceed, we need to discuss certain elements of the modular group. Recall the generators~$\gv$ and~$\gt$ of the modular group from Paragraph~\ref{CentCharge}. 
We now introduce the following matrices (cf.~\cite{Frasch}, \S~1, Eq.~(6), p.~230), which implicitly appear already in Proposition~\ref{RelAct}:
\begin{defn}
For $q,r \in \Z$, we define
$\gd(q,r):= \gv \gt^{r} \gv^{-1} \gt^q \gv \gt^{r}$.
\end{defn}
Explicitly, these matrices are
\begin{align*}
\gd(q,r) &= 
\begin{pmatrix} 0 & -1 \\ 1 & 0 \end{pmatrix} 
\begin{pmatrix} 1 & r \\ 0 & 1 \end{pmatrix}
\begin{pmatrix} 0 & 1 \\ -1 & 0 \end{pmatrix} 
\begin{pmatrix} 1 & q \\ 0 & 1 \end{pmatrix}
\begin{pmatrix} 0 & -1 \\ 1 & 0 \end{pmatrix} 
\begin{pmatrix} 1 & r \\ 0 & 1 \end{pmatrix} \\
&= 
\begin{pmatrix} 0 & -1 \\  1 & r \end{pmatrix}
\begin{pmatrix} 0 & 1 \\ -1 & -q\end{pmatrix}
\begin{pmatrix} 0 & -1  \\ 1 & r\end{pmatrix} \\
&= 
\begin{pmatrix} 1 & q \\ -r & 1-qr \end{pmatrix}
\begin{pmatrix} 0 & -1  \\ 1 & r\end{pmatrix}
= \begin{pmatrix} q & qr-1 \\ 1-qr & r(2-qr) \end{pmatrix}
\end{align*}
From this we see that, if $qr \equiv 1 \pmod{N}$, then
the image of $\gd(q,r)$ in~$\SL(2,\Z_N)$ under reduction modulo~$N$
becomes diagonal. In any case, we will
denote the image of~$\gd(q,r)$ under reduction modulo~$N$ 
by~$\bar{\gd}(q,r)$.

Denoting the transpose of a $2\times 2$-matrix~$A$ by~$A^T$, we have the following lemma, which appears in \cite{Beyl}, Sec.~1, no.~(1.5), p.~24:
\begin{lemma}
For all $g \in \SL(2,\Z)$, we have $\gv g^{-1} = g^T \gv$. 
\end{lemma}
\begin{pf}
The assertion can be reformulated by saying that the
two group automorphisms $g \mapsto \gv g \gv^{-1}$ 
and $g \mapsto g^{-1T}$ are equal, a fact that can be verified on generators. For $g=\gv$, this is obvious, and for $g=\gt$, the corresponding calculation can be found in \cite{SoZhu}, Par.~1.1, p.~10. 
It is also easy to check the assertion directly.
\qed
\end{pf}

For the matrices $\gd(q,r)$, this lemma implies the following commutation relations:
\begin{prop}
$\gv \gd(q,r)^{-1} = \gd(-q,-r) \gv^{-1} = \gd(q,r)^T \gv$
\end{prop}
\begin{pf}
For the first identity, we have
\begin{align*}
\gv \gd(q,r)^{-1} 
= \gv \gt^{-r} \gv^{-1} \gt^{-q} \gv \gt^{-r} \gv^{-1}
= \gd(-q,-r) \gv^{-1}
\end{align*}
The second identity follows from the preceding lemma.
\qed
\end{pf}

Because the transpose of a diagonal matrix is equal to itself, 
this proposition implies in particular that
$\gv \gd(q,r)^{-1} \equiv \gd(q,r) \gv \pmod{\Gamma(N)}$ if
$qr \equiv 1 \pmod{N}$.

\subsection[Application of Beyl's theorem]{} \label{BeylAppl}
We now consider a projective congruence datum over an algebraically closed field~$K$ of characteristic zero. As explained in Paragraph~\ref{DefModData}, we denote its exponent by~$N$ and its normalized exponent by~$N_o$. We then have 
a group homomorphism $\alpha: \SL(2,\Z_{N_o}) \rightarrow \PGL(m,K)$ 
that maps the reduced matrices~$\bar{\gv}$ and~$\bar{\gt}$ to the cosets of~$\V$ resp.~$\T$. As we saw in Paragraph~\ref{GroupCohom}, we can lift~$\alpha$ to a group homomorphism $\beta: \SL(2,\Z_{N_o}) \rightarrow \Pn_2 \GL(m,K)$. We choose two representatives $\V'', \T'' \in \GL(m,K)$ for the images of the generators, so that we have
$$\beta(\bv) = \{\pm \V''\} \qquad \qquad \beta(\bt) = \{\pm \T''\}$$
Because~$\beta$ lifts~$\alpha$, there are 
nonzero numbers~$s,t \in K^\times$ such that
$$\V'' = \frac{1}{s} \V \qquad \qquad \T'' = \frac{1}{t} \T$$
Applying~$\beta$ to the defining relations of the modular group described in Paragraph~\ref{CentCharge}, we get $\beta(\bv)^4 = 1$ and $(\beta(\bt) \beta(\bv))^3 = \beta(\bv)^2$. Therefore, there are signs 
$v_s, v_t \in \{1,-1\}$ such that
$$\V''^4 = v_s \E \qquad (\T'' \V'')^3 = v_t \V''^2$$
In terms of the parameters~$s$ and~$t$, this equation can be
written as
$$\V^4 = v_s s^4 \E \qquad (\T \V)^3 = v_t s t^3 \V^2$$
In Paragraph~\ref{CentCharge}, we have rewritten part of the definition of
a modular datum in the form
$\V^2 = n \CC$ resp.\ $\frac{n_o}{t_o^2} (\T \V)^3 = \gp  \V^2$.
Comparing this with the preceding equations, we get
$$v_s s^4 = n^2 \qquad \qquad v_t s t^3 = \frac{t_o^2}{n_o} \gp$$
From this, we see that the new parameters~$s$ and~$t$ are essentially roots of unity:
\begin{prop}
For an extended projective congruence datum that is also Galois, we have that $s/D$ is an eighth root of unity and that~$t$ is a 72th root of unity, where~$D$ is the generalized rank of the extended modular datum. Moreover, for the central charge~$\ell$ of the datum, we have
$$t^3 = v_t \frac{D}{s} t_o^3 \ell^3$$
\end{prop}
\begin{pf}
From the first equation above, we get $s^4 = v_s n^2 = v_s D^4$, so that
$(s/D)^4 = v_s$ and thus $(s/D)^8 = 1$. Dividing the second equation above by~$D$, we get 
$$v_t \frac{s}{D} t^3 = \frac{t_o^2}{n_o} \frac{\gp}{D} = t_o^3 \ell^3$$
where the second equation uses Definition~\ref{CentCharge}. Squaring this and using the proof of Proposition~\ref{CentCharge}, we get
$$(\frac{s}{D})^2 t^6 = \pm t_o^4 \frac{\gp}{\gm}$$
Now $\gp/\gm$ is a twelfth root of unity by Proposition~\ref{12root},
and we have also discussed there that $t_o^{24}=1$. Combining all this, we get 
$$t^{72} = (\pm t_o^4 \frac{\gp}{\gm})^{12} = 1$$
as asserted.
\qed
\end{pf}

From this proposition, we see that the matrices~$\V''$ and~$\T''$ differ from the homogeneous matrices~$\V'$ and~$\T'$ introduced in Definition~\ref{CentCharge} only by roots of unity: $\V'$ and~$\V''$ differ by the eighth root of unity~$s/D$, and $\T'$ and~$\T''$ differ by the 24th root of unity~$t/(t_o \ell)$. While the matrices~$\V'$ and~$\T'$
lead to an ordinary representation of the modular group, the matrices~$\V''$ and~$\T''$ only lead to a representation modulo signs,
but this latter homomorphism to~$\Pn_2\GL(m,K)$ has the advantage to factor over the reduced modular group~$\SL(2,\Z_{N_o})$.

\subsection[Fourth roots of unity]{} \label{4root}
In the last paragraph, we have associated with a projective congruence datum over an algebraically closed field~$K$ of characteristic zero the group homomorphism~$\beta: \SL(2,\Z_{N_o}) \rightarrow \Pn_2\GL(m,K)$, the matrices~$\V''$ and~$\T''$, and the numbers~$s$ and~$t$. For these quantities, we have the following analogue of
Proposition~\ref{RelAct}:
\begin{prop}
Consider an extended projective congruence datum that is also Galois,
and suppose that we are given two integers $q,q' \in \Z$ satisfying $qq' \equiv 1 \pmod{72N}$. Then we have
$$\beta(\bar{\gd}(q,q')) = 
\{ \pm (\frac{s}{D})^{q-1} t^{2(q-q')} \; P(\sigma_q^{-1})\}$$
where~$D$ is the generalized rank of the extended modular datum.
\end{prop}
\begin{pf}
Dividing the equation in Proposition~\ref{RelAct} by $t^{2q'+q} s$, we get
\begin{align*}
\V'' \T''^{q'} \V''^{-1} \T''^q \V'' \T''^{q'} &=  
\frac{t_o^{2q}}{n_o} \frac{\sigma_q(\gp)}{t^{2q'+q} s} \; P(\sigma_q^{-1}) 
\end{align*}
A slight modification of the argument given in the proof of Proposition~\ref{12root} shows that $D \in \Q_{4N}$, so that $\sigma_q(D)$
is well-defined. Now determine $v \in \{1,-1\}$ such that $\sigma_q(D) = vD$. Using  Proposition~\ref{BeylAppl} and the formulas stated before it, we get
$$\frac{t_o^{2q}}{n_o} \sigma_q(\gp) = \sigma_q(\frac{t_o^{2}}{n_o} \gp) = \sigma_q(v_t D \frac{s}{D} t^3) = 
v_t v D (\frac{s}{D})^q t^{3q}$$
Inserting this into the preceding equation, it becomes
$$\V'' \T''^{q'} \V''^{-1} \T''^q \V'' \T''^{q'} =  
v_t v (\frac{s}{D})^{q-1} t^{2(q-q')} \; P(\sigma_q^{-1})$$
But according to Definition~\ref{DiagMat}, this is what we have
to verify.
\qed
\end{pf}

This enables us to take the fourth step in the sharpening of Proposition~\ref{FusSym}:
\begin{thm}
For an extended projective congruence datum that is also Galois, we have $\ell^{24} = 1$ and $\gp^4 = t_o^8 \gm^4$.
\end{thm}
\begin{pf}
As these identities clearly still hold if we enlarge the base field,
we can assume that~$K$ is algebraically closed, so that the preceding
discussion applies. Recall from the proof of Proposition~\ref{CentCharge} that $t_o^2 \ell^6 = \pm \gp/\gm$, so that our two assertions are equivalent.
Moreover, we can assume that~$N$ is even, because otherwise
we have already~$\gp = \pm t_o^2 \gm$ by Theorem~\ref{OddExp}.

We claim that we can also assume that~$N$ is divisible by~$3$:
From Proposition~\ref{12root}, we know that 
$(\gm/\gp)^{12} = \fs(-1)^{12} = 1$, and we also discussed there that $t_o^{24}=1$. Furthermore, we have 
$(\gm/\gp)^{N} =1$ by Proposition~\ref{GaussSum},
and $t_o^N=1$ by definition. If~$N$ is not divisible by~$3$,
then these equations imply together that 
$(\gm/\gp)^{4} = 1$ and $t_o^8=1$, in which case our assertions hold.

In the remaining cases, $N$ and~$72N$ therefore have the same prime divisors, so that, if~$q$ is relatively prime to~$N$, we can choose $q' \in \Z$ such that $qq' \equiv 1 \pmod{72N}$. Since~$N_o$ divides~$N$, we then also have $qq' \equiv 1 \pmod{N_o}$, so that by Proposition~\ref{DiagMat} the matrices $\gd(-q,-q')$ and $\gd(q,q') \gv^2$ reduce in~$\SL(2,\Z_{N_o})$ to the same diagonal matrix. We have 
$\beta(\bar{\gv}) = \{\pm \V''\} = \{\pm \frac{1}{s} \V \}$,
therefore
\mbox{$\beta(\bar{\gv}^2) = \{\pm \frac{1}{s^2} \V^2 \} = \{\pm \frac{n}{s^2} \CC \}$},
and so the preceding proposition yields
\begin{align*}
\{\pm (\frac{s}{D})^{-q-1} t^{2(q'-q)} \; P(\sigma_{-q}^{-1})\}
&= \beta(\bar{\gd}(-q,-q')) \\
&= \beta(\bar{\gd}(q,q') \bar{\gv}^2)
= 
\{ \pm (\frac{s}{D})^{q-1} t^{2(q-q')} \frac{n}{s^2}\; P(\sigma_q^{-1}) \CC\}
\end{align*}
Since $P(\sigma_{-1}) = P(\gamma) = \CC$, we have $P(\sigma_{-q}^{-1})=P(\sigma_q^{-1}) \CC$. As we have $D^2=\pm n$,
this implies that there is a sign $v \in \{1,-1\}$ such that
\begin{align*}
(\frac{s}{D})^{-q-1} t^{2(q'-q)}  &= 
v (\frac{s}{D})^{q-1} t^{2(q-q')} \frac{D^2}{s^2} 
= v (\frac{s}{D})^{q-3} t^{2(q-q')} 
\end{align*}
or alternatively $t^{4(q-q')} = v(\frac{s}{D})^{-2q+2}$.
Now we know from Proposition~\ref{BeylAppl} that~$s/D$ is an eighth root of unity, so that by raising the last equation to the fourth power we get
$t^{16(q-q')} = 1$. The same proposition yields that~$t^{72}=1$, which
means that already~$t^{8(q-q')} = 1$. Raising this equation to the~$q$-th power, it becomes~$t^{8(q^2-1)} = 1$, or~$t^{8q^2} = t^8$.

All of this shows that~$t^8$ is a $72N$-th root of unity which is invariant
under~$\sigma_q^2$ for every~$q$ that is relatively prime to~$72N$. So the `definition of~24', given in Lemma~\ref{12root}, implies that 
$(t^8)^{24}=1$. As we also have $t^{72}=(t^8)^{9}=1$, and $3=\gcd(24,9)$, we get $t^{24}=(t^8)^{3}=1$. Using Proposition~\ref{BeylAppl} once again, we see that
$t^3 = v_t \frac{D}{s} t_o^3 \ell^3$ is an eighth root of unity, and as
$\frac{s}{D}$ and $t_o^3$ are also eighth roots of unity, $\ell^3$ must be an eighth root of unity, too, which is the assertion.
\qed
\end{pf}

If the modular datum is in addition normalized, the theorem obviously
yields that $\gp^4=\gm^4$. Furthermore, we have explained in Paragraph~\ref{CentCharge} that in the case $K=\C$ we can write $\ell = e^{2\pi i c/24}$ by using the additive central charge~$c$. The equation $\ell^{24}=1$ then means that~$c$ must be an integer. 

We note that in~\cite{CostGann}, A.~Coste and T.~Gannon deduce the same result from a slightly different hypothesis (cf.~\cite{CostGann}, \S~2.4, Prop.~3.b, p.~9). To understand the relation of the two results, let us consider a normalized modular datum. We can extend it by choosing a generalized rank~$D$ and a multiplicative central charge~$\ell$, and then form the homogeneous matrices~$\V'=\V/D$ and~$\T'=\T/\ell$. Following their treatment in this situation, we would then assume that~$\T'$, and not~$\T$, satisfies the Galois condition. As the $(o,o)$-component of~$\T'$
is~$1/\ell$, the fact that~$\ell^{24}=1$ then follows as in our remark
after Lemma~\ref{12root}; in fact, our argument there was directly adapted
from the argument of Coste and Gannon. However, if we impose the Galois condition on~$\T$ instead of~$\T'$, then~$\T'$ will also satisfy the Galois condition if and only if~$\ell^{24}=1$, so that the argument given by Coste and Gannon becomes unavailable. But for the applications to semisimple Hopf algebras that we have in mind, it is the Galois condition on~$\T$ that is needed.

We conclude this paragraph with a brief table that 
summarizes the main results that we have obtained in short form:
\begin{center}
\setlength{\extrarowheight}{2.3pt}
\begin{tabular}[t]{|c|c|c|}\hline 
Result & Assumptions & Paragraph \\ \hline 
$\gp^{2N}=\gm^{2N}$, $t_o^{N}=1$ & integral modular datum & \ref{GaussSum} (cf.~\ref{FusSym}) \\ \hline 
$\gp^{12}=\gm^{12}$, $t_o^{24}=1$ & Galois modular datum & \ref{12root} \\ \hline 
$\gp^{2} = t_o^4 \gm^{2}$ & integral modular datum, $N$ odd & \ref{OddExp}  \\ \hline 
$\gp = (-1)^{\frac{n-1}{2}} t_o^2 \gm$ & integral modular datum, $N$ and $n$  odd & \ref{ClassGauss}  \\ \hline 
$\gp^{4} = t_o^8 \gm^{4}$ & Galois projective congruence datum  & \ref{4root} \\ \hline 
\end{tabular} 
\end{center}

\subsection[Cauchy's theorem]{} \label{Cauchy}
As we already mentioned there, the hypothesis of Theorem~\ref{ClassGauss} raises the question how the prime factors of~$n$ and~$N$ are related. As we will explain in Paragraph~\ref{FactHopf}, every semisimple
factorizable Hopf algebra leads to a modular datum whose global
dimension is just the dimension of the Hopf algebra. In this situation,
Cauchy's theorem for Hopf algebras (cf.~\cite{YYY2}, Thm.~3.4, p.~26) asserts that a prime that divides~$n$ also divides~$N$. For integral modular data, we have the following weak version of this result:
\begin{prop}
For an integral modular datum, suppose that an odd prime~$p$ divides~$n$ an odd number of times. Then it also divides~$N$.
\end{prop}
\begin{pf}
Let~$l$ be the product of all primes that appear an odd number of times in the prime factorization of~$n$. Then~$l$ is squarefree, and in the
factorization of~$n/l$, each prime appears an even number of times, so
that we can write it as a complete square $n/l=k^2$ for some positive
integer~$k$. Our assumption now means that~$p$ divides~$l$, so that 
in particular~$l \neq 1$.

By possibly enlarging the base field, we can choose a rank for the modular datum, i.e., a number $D \in K$ such that~$D^2=n$. Recall from the first step in the proof of Proposition~\ref{12root} that $D \in \Q_{4N}$, so that
the cyclotomic field $\Q_{4N}$ contains the quadratic number field $\Q(D)$,
which is clearly isomorphic to the quadratic number field $\Q(\sqrt{n}) = \Q(\sqrt{l}) \subset \C$ and therefore has discriminant~$4l$ if
$l \equiv 2 \pmod{4}$ or $l \equiv 3 \pmod{4}$, and discriminant~$l$ if
$l \equiv 1 \pmod{4}$ (cf.~\cite{N}, Kap.~I, \S~2, Aufg.~4, p.~16). In any case, its discriminant is divisible by~$p$.  

By the discriminant tower theorem (cf.~\cite{N}, Kap.~III, \S~2, Kor.~(2.10), p.~213), the discriminant of~$\Q(D)$ divides the discriminant of~$\Q_{4N}$. But a prime that divides the discriminant of~$\Q_{4N}$ also divides~$4N$ (cf.~\cite{Wash}, Chap.~2, Prop.~2.7, p.~12; see also \cite{N}, Kap.~I, \S~10, Lem.~(10.1), p.~62). Therefore, $p$ divides~$4N$, and therefore~$N$.
\qed
\end{pf}

It is also possible to say something about the prime~$p=2$:
\begin{corollary}
Consider a projective congruence datum that is also Galois. If~$N$ denotes its exponent and~$n$ its global dimension,
then we have $N \equiv 0 \pmod{4}$ if $n \equiv 2 \pmod{4}$.
\end{corollary}
\begin{pf}
We can assume that the base field~$K$ is algebraically closed, and therefore extend the datum by choosing a generalized rank~$D$ and a 
central charge~$\ell$. We then know from Theorem~\ref{4root} that
$t_o^8 \gm^4 = \gp^4$, and therefore $t_o^4 \gm^2 = \pm \gp^2$. If 
$t_o^4 \gm^2 = - \gp^2$,
then~$\Q_N$ contains the primitive fourth root of unity~$t_o^2 \gm/\gp$. If now~$N$ were odd, we would have 
$1 = (t_o^2 \gm/\gp)^{2N} = (-1)^N = -1$
by~\cite{Wash}, Exerc.~2.3, p.~17, so~$N$ must be even. But this implies the stronger equation $(t_o^2 \gm/\gp)^{N} = 1$, which shows that~$N$ is divisible by~$4$.

We can therefore assume that $t_o^4 \gm^2 = \gp^2$, so that 
$t_o^2 \gm = v \gp$ for a sign $v \in \{1,-1\}$. According to Proposition~\ref{DefModData}.\ref{DefModDataRecGauss}, this means that
$\gp^2 = v t_o^2 \gp \gm = v t_o^2 n n_o^2$. As in the proof of the preceding proposition,
we write~$n = k^2 l$, where~$l$ is squarefree. The last equation then implies that $(\gp/(t_o n_o k))^2 = v l$.

Because $n \equiv 2 \pmod{4}$, we must have $l = 2 l'$ for an odd integer~$l'$, and~$l'$ divides~$N$ by the preceding proposition.
If we denote the classical Gaussian sum coming from a primitive $l'$-th
root of unity by~$\Gp$, then we know from Lemma~\ref{ClassGauss}
that $\Gp^2 = (-1)^{\frac{l'-1}{2}} l'$; moreover, we have
$\Gp \in \Q_{l'} \subset \Q_N$. Combining this with the previous equation, we get
$$(\frac{\gp}{t_o n_o k \Gp})^2 = (-1)^{\frac{l'-1}{2}} \frac{vl}{l'} = \pm 2$$
On the other hand, consider a primitive eighth root of unity~$\xi \in K$.
Then we have $\xi^4 = -1$ and therefore
$$(\frac{1+\xi^2}{\xi})^2 = 2 \qquad \text{and} \qquad (\frac{1-\xi^2}{\xi})^2 = -2$$
This shows that 
$$\frac{\gp}{t_o n_o k \Gp} = \pm \frac{1+\xi^2}{\xi} \qquad \text{or} \qquad \frac{\gp}{t_o n_o k \Gp} = \pm \frac{1-\xi^2}{\xi}$$
In any case, we have $\gp/(t_o n_o k \Gp) \in \Q_N \cap \Q_8$. If~$N$ is odd,
then $\Q_N \cap \Q_8 = \Q$ by \cite{Wash}, Prop.~2.4, p.~11. If
$N \equiv 2 \pmod{4}$, we have $N=2N'$ for an odd number~$N'$,
so that $\Q_N = \Q_{N'}$ and again $\Q_N \cap \Q_8 = \Q$ by the same argument. In both cases, we get that $\gp/(t_o n_o k \Gp)$ is a rational number whose square is~$\pm 2$, which is a contradiction. Therefore,
$N \equiv 0 \pmod{4}$.
\qed
\end{pf}

Although the preceding argument is sufficient to prove the assertion, it is worth noting that in the second case, where $\gp^2 = t_o^4 \gm^2$, we actually get the stronger statement that $N \equiv 0 \pmod{8}$. To see this, note that, because we now know that~$N$ is divisible by~$4$, we have
$\xi^2 \in \Q_N$, and we have seen in the proof that
$$\frac{1+\xi^2}{\xi} \in \Q_N \qquad \text{or} \qquad \frac{1-\xi^2}{\xi} \in \Q_N $$
Both equations clearly imply that $\xi \in \Q_N$, so that~$N$ must be divisible by~$8$.

This fact is worth pointing out because it appears that, under the assumptions made in the corollary, we always have that~$N$ divides~$n^2$.
We are therefore lead to conjecture the following stronger form of the corollary:
\begin{conj}
Consider a projective congruence datum that is also Galois. If~$N$ denotes its exponent and~$n$ its global dimension, then we have 
\begin{enumerate}
\item 
$N \equiv 4 \pmod{8}$

\item
$\gp^2 = -t_o^4 \gm^2$
\end{enumerate}
if $n \equiv 2 \pmod{4}$.
\end{conj}

The preceding argument at least shows that the first statement implies the second. Note that in the normalized case the second statement asserts that
$\gm/\gp$ is a primitive fourth root of unity. We will see in Paragraph~\ref{Semion} an example where all the assumptions, and therefore all the conclusions, made in this paragraph are satisfied: We there have that $n=2$, $N=4$, $t_o=1$, and $\gp^2 = - \gm^2$.

\newpage
\section{Hopf Algebras} \label{Sec:HopfAlg}
\subsection[Modular categories]{} \label{ModCat}
Modular data axiomatize certain properties of modular categories,
as we explain now. A modular category in the sense of V.~G.~Turaev (cf.~\cite{Tur}, Sec.~II.1.4, p.~74) comes with a finite family 
$(V_i)_{i \in I}$ of objects such that for one index $o \in I$ the corresponding object~$V_o$ is the unit object. It follows from the axioms of a modular category that $K:=\End(V_o)$ is always a commutative ring (cf.~\cite{Tur}, Sec.~II.1.1, p.~72); we assume here that it is an algebraically closed field of characteristic zero. The duality in the category leads to an involution on the index set with the property that $o^*=o$ (cf.~\cite{Tur}, p.~75). If the quasisymmetry of the category is denoted by~$c$, then we define the Verlinde matrix~$\V$ via
$$s_{ij} = \Tr_q(c_{V_{j^*},V_i} \circ c_{V_i,V_{j^*}})$$
where $\Tr_q$ denotes the categorical trace, or quantum trace, as it is sometimes called (cf.~\cite{Tur}, Sec.~I.1.5, p.~21; \cite{Kas}, Def.~XIV.4.1, p.~354). Note that here we have deviated slightly from the convention in \cite{Tur}, p.~74; we include in Paragraph~\ref{FactHopf} below a table that compares the conventions of a couple of references. The Dehn matrix~$\T$ arises from the twist of the category, which acts on the simple object~$V_i$ via multiplication by a scalar~$t_i \in K$; this number is denoted by~$v_i$ in \cite{Tur}, p.~76.

This explains how the structure elements of a modular datum arise from a modular category; however, the axioms that we have given in Definition~\ref{DefModData} will only be satisfied if the category is in addition semisimple. By this, we mean that the category is \mbox{$K$-linear}, that the objects~$V_i$ above are simple, and that every object is isomorphic to a finite direct sum of some of the objects~$V_i$, possibly with repetitions (cf.~\cite{Tur}, Sec.~II.4.1, p.~99). For these categories, the fact that the numbers~$t_i$ are roots of unity is known as Vafa's theorem (cf.~\cite{BakKir}, Thm.~3.1.19, p.~57). The equation $t_{i^*}=t_i$ holds by \cite{Tur}, Sec.~II.3.3, Eq.~(3.3.b), p.~90.
As noted in \cite{Tur}, Sec.~II.1.4, p.~74f, the Verlinde matrix is always symmetric, and~$s_{io}$ is always equal to the categorical dimension of the $i$-th simple object, which is shown there to be nonzero. In view of \cite{Tur}, Exerc.~II.1.9.2, p.~78, all of this also holds after our modification.

Our Axiom~\ref{DefModData}.\ref{DefModData3} about the form of the inverse of~$\V$ is in this setting satisfied by~\cite{Tur}, Sec.~II.3.8, Eq.~(3.8.a), p.~97, resp.  Sec.~II.3.9, p.~98. Similarly, Axiom~\ref{DefModData}.\ref{DefModData4} is satisfied in this situation by~\cite{Tur}, Sec.~II.3.8, Eq.~(3.8.c), p.~97. Turaev's version is in fact closer to our constant form of Axiom~\ref{DefModData}.\ref{DefModData4}. Although he gives no name to the Gaussian sum, he also introduces an element that corresponds to our reciprocal Gaussian sum~$\gm$ and denotes it by~$\Delta$ (cf.~\cite{Tur}, Sec.~II.1.6, p.~76). However, the term `Gaussian sum' is often used for this element (cf.~\cite{Mue2}, Def.~1.1, p.~160), as it reduces to the classical Gaussian sum in certain examples that we will consider in Paragraph~\ref{ExamRadf}. This is also pointed out by Turaev when he
derives his variant of Proposition~\ref{DefModData}.\ref{DefModDataRecGauss}
(cf.~\cite{Tur}, Sec.~II.2.4, Eq.~(2.4.a), p.~83). Finally, Axiom~\ref{DefModData}.\ref{DefModData5} about the numbers~$N_{ij}^k$ is a property of modular categories known as the Verlinde formula
(cf.~\cite{Tur}, Thm.~II.4.5.2, p.~106; \cite{BakKir}, Thm.~3.1.14, p.~54), which shows that this number is the multiplicity of~$V_k$ in the decomposition of~$V_i \o V_j$ into simple objects, and therefore in particular a nonnegative integer.

Modular data coming from modular categories are always normalized:
The fact that $t_o=1$ is shown in \cite{Tur}, Sec.~I.1.2, p.~20, and the equation $n_o=1$ is shown in \cite{Tur}, Lem.~I.1.5.1, p.~22. 
However, they are not always integral. We therefore now proceed to
study a class of semisimple modular categories that lead to integral modular data, namely those coming from factorizable Hopf algebras.

\subsection[Factorizable Hopf algebras]{} \label{FactHopf}
So, we now consider a semisimple Hopf algebra~$A$ with coproduct~$\da$,
counit~$\ea$, and antipode~$\sa$ over an algebraically closed field~$K$ of characteristic zero. We denote the dimension of~$A$ by~$n$; note that semisimple Hopf algebras are always finite-di\-mensional (cf.\ \cite{Sw1}, Cor.\ 2.7, p.~330; \cite{Sw}, Chap.\ V, Exerc.\ 4, p.\ 108). Furthermore, $A$ is also cosemisimple and the antipode is an involution (cf.~\cite{LR1}, Thm.~4, p.~195; \cite{LR2}, Thm.~3.3, p.~276).

We also assume that~$A$ is quasitriangular. This means that $A$ possesses a so-called R-matrix, i.e., an invertible element $R=\sum_{i=1}^k a_i \o b_i \in A \o A$ that satisfies
$\da^{\cop}(a) = R \da(a) R^{-1}$ as well as
$$(\da \o \id)(R) = \sum_{i,j=1}^k a_i \o a_j \o b_i b_j \qquad
(\id \o \da)(R) = \sum_{i,j=1}^k a_i a_j \o b_j \o b_i$$
(cf.~\cite{Kas}, Def.~VIII.2.2, p.~173; \cite{M}, Def.~10.1.5, p.~180; \cite{Tur}, Sec.~XI.2.1, p.~496). From the R-matrix, we derive two further elements, namely
$$R':=\sum_{i=1}^k b_i \o a_i \qquad  \qquad 
\ua:=\sum_{i=1}^k \sa(b_i) a_i$$
The element~$\ua$ is called the Drinfel'd element; it is always invertible.
$A$~is called factorizable if the tensor~$R'R$ has maximal rank, i.e., if it cannot be written as a sum of less than~$n$ decomposable tensors.

The module category of a semisimple factorizable Hopf algebra is modular:
For the objects $(V_i)_{i \in I}$, the Wedderburn structure theorem provides a system of representatives for the isomorphism classes of simple modules with the required properties. The object~$V_o$ is chosen to be the base field~$K$, turned into an $A$-module via the counit. As explained in \cite{Tur}, the duality on the category arises from the antipode~$\sa$ (cf.~\cite{Tur}, Sec.~XI.1.3, p.~494f; see also \cite{Kas}, Sec.~XIV.2, Ex.~1, p.~347), and the quasisymmetry of the category arises from the R-matrix (cf.~\cite{Tur}, Sec.~XI.2.3, p.~498f; see also \cite{Kas}, Prop.~XIII.1.4, p.~318).

Because in our situation the antipode is an involution, the inverse Drinfel'd element~$\ua^{-1}$ can be used as a ribbon element. For this ribbon element, the categorical trace coincides with the usual trace (cf.~\cite{Tur}, Lem.~XI.3.3, p.~501; see also \cite{Kas}, Prop.~XIV.6.4, p.~363), so that the categorical dimensions coincide with the ordinary dimensions. This implies that the entries of the Verlinde matrix~$\V$ are given as
$$s_{ij} = (\chi_i \o \chi_{j^*})(R'R) $$
where $\chi_i$ denotes the character of~$V_i$. The diagonal entry~$t_i$
of the Dehn matrix~$\T$ is determined by the condition that~$\ua^{-1}$
acts on~$V_i$ as $t_i \id_{V_i}$. The decisive axiom of a modular category to check now is the invertibility of the Verlinde matrix; this axiom follows
from factorizability by establishing a version of Axiom~\ref{DefModData}.\ref{DefModData3} (cf.~\cite{SchneiderFact}, Rem.~3.4, p.~1895; \cite{SoZhu}, Prop.~5.3, p.~48).

Because the module category of a semisimple factorizable Hopf algebra is modular, the discussion in Paragraph~\ref{ModCat} implies that
the quintuple $(I,o,{}^*,\V,\T)$ is a normalized modular datum. Of course,
it is also possible to deduce this directly from the definition of a factorizable Hopf algebra (cf.~\cite{SoZhu}, Sec.~5). But the important
new feature now is that these modular data are always integral: Because
$\chi_o = \ea$, it follows from fundamental properties of R-matrices (cf.~\cite{M}, Prop.~10.1.8, p.~180; \cite{Kas}, Thm.~VIII.2.4, p.~175; \cite{Tur}, Lem.~XI.2.1.1, p.~497) that
$$s_{io} = (\chi_i \o \ea)(R'R) = \chi_i(\A) = \dim(V_i)$$
so that $n_i = s_{io} = \dim(V_i)$ is the dimension of a vector space,
and therefore a positive integer. 

The various quantities that we have associated with a modular datum reappear in a slightly different form in the Hopf algebra situation.
It follows directly from the Wedderburn structure theorem that the global dimension~$n$ of the modular datum is equal to the dimension of the Hopf algebra. The exponent~$N$ of the modular datum is obviously equal to the order of the Drinfel'd element~$\ua$. It is also equal to the exponent as originally defined in terms of Sweedler powers (cf.~\cite{KashAnti}, p.~1261; \cite{KashPow}, p.~159; \cite{EG4}, Def.~2.1, p.~132), although this is not 
completely obvious: It follows by combining \cite{EG4}, Thm.~2.5, p.~133; \cite{KashPow}, Thm.~3.4, p.~170, and \cite{SoZhu}, Lem.~3.1, p.~22.

The Gaussian sum $\gp = \sum_{i \in I} n_i^2 t_i$ also admits another description in the Hopf algebra setting: It again follows from the Wedderburn structure theorem that, for an element $a \in A$, the trace~$\chi_R(a)$ of the left multiplication by~$a$ is equal 
to~$\chi_R(a) = \sum_{i \in I} n_i \chi_i(a)$.
The function~$\chi_R$ is called the character of the regular representation;
it is also a left and right integral on~$A$ (cf.~\cite{LR2}, Prop.~2.4, p.~273). Because $\ua^{-1}$ acts on~$V_i$ as $t_i \id_{V_i}$, we have
$\chi_i(\ua^{-1}) = n_i t_i$, which implies that $\gp = \chi_R(\ua^{-1})$; similarly, we see that $\gm = \chi_R(\ua)$. As a consequence, the fusion symbol~$\fs(q)$ reduces to the Hopf symbol~$\jac{q}{A}$ in this case (cf.~\cite{SoZhu}, Def.~12.1, p.~114),
and Proposition~\ref{DefModData}.\ref{DefModDataRecGauss} becomes
$\chi_R(\ua^{-1}) \chi_R(\ua)= \dim(A)$ (cf.~\cite{SoZhu}, Par.~5.3, p.~49). 

To facilitate the comparison between the different conventions used in
the various references, we include the following table:
\begin{center}
\setlength{\extrarowheight}{1.5pt}
\begin{tabular}{|c|c|c|c|}
\hline 
Present notation & \cite{Tur} & \cite{BakKir} & \cite{SoZhu}\\ 
\hline 
$s_{ij}$ & $S_{ij^*}$ (p.~74) & $\tilde{s}_{ij}$  (p.~47) & $s_{ij^*}$  (p.~45)\\ 
\hline 
$t_i$ & $v_i$ (p.~76)  & $\theta_i$ (p.~44) & $1/u_i$ (p.~45)\\ 
\hline 
$n_i$ & $\dim(i)$ (p.~74)  & $d_i$ (p.~44) & $n_i$ (p.~43) \\ 
\hline 
$\gp$ & $\Delta_{\bar{\cal V}}$ (p.~21)  & $p_+$ (p.~49) & $\chi_R(\ua^{-1})$ (p.~43) \\ 
\hline
$\gm$ & $\Delta = \Delta_{\cal V}$ (p.~77)  & $p_-$ (p.~49) & $\chi_R(\ua)$ (p.~43) \\ 
\hline 
$D$ & $\cal D$ (p.~76)  & $D$ (p.~51) & $1/\kappa$ (p.~90) \\ 
\hline
$\ell$ &  & $\zeta$ (p.~51) &  \\ 
\hline 
$N_{ij}^k$ & $h^{ij}_k$ (p.~105)  & $N_{ij}^k$ (p.~44) &  \\
\hline 
\end{tabular} 
\end{center}

\subsection[The main result]{} \label{MainResult}
We are now in the position to present the main result of this article: 
\begin{thm}
Consider a semisimple factorizable Hopf algebra~$A$ over a field~$K$ of characteristic zero. Let~$n$ be the dimension of~$A$, $\ua$ its Drinfel'd element, and~$\chi_R$ the character of its regular representation.
Then the following holds:
\begin{enumerate}
\item 
If $n$ is odd, then we have
$$\chi_R(u^{-1}) = 
\begin{cases}
\chi_R(u) &: n \equiv 1 \pmod{4} \\
- \chi_R(u) &: n \equiv 3 \pmod{4}
\end{cases}$$
Moreover, the Hopf symbol coincides with the Jacobi symbol; i.e.,
we have
$$\jac{q}{A} = \jac{q}{n}$$
for all~$q \in \Z$.

\item 
If $n$ is even, then we have $\chi_R(u^{-1})^4 = \chi_R(u)^4$.
\end{enumerate}
\end{thm}
\begin{pf}
Because these identities will still hold when we enlarge the base field, we can assume that $K$ is algebraically closed, so that the discussion in Paragraph~\ref{FactHopf} applies. We know that the exponent~$N$ of a semisimple Hopf algebra is odd if and only if its dimension~$n$ is odd (cf.~\cite{EG4}, Thm.~4.3, p.~136; \cite{YYY1}, Cor.~4, p.~93), and therefore the equation $\chi_R(u^{-1}) = (-1)^{(n-1)/2} \chi_R(u)$
in the first assertion follows directly from Theorem~\ref{ClassGauss}. Even stronger, Cauchy's theorem for Hopf algebras implies that~$N$ and~$n$ have the same prime divisors (cf.~\cite{YYY2}, Thm.~3.4, p.~26), so that an integer~$q$ is relatively prime to~$N$ if and only if it is relatively prime to~$n$.
Therefore, the equation $\jac{q}{A} = \jac{q}{n}$ also follows from Theorem~\ref{ClassGauss}.

The second assertion follows similarly from Theorem~\ref{4root}: By the projective congruence subgroup theorem, modular data coming from semisimple factorizable Hopf algebras are projective congruence data, and they are also Galois (cf.~\cite{SoZhu}, Thm.~9.4, p.~94, Lem.~12.2, p.~115).
\qed
\end{pf}

We have already discussed in Paragraphs~\ref{OddExp},~\ref{ClassGauss}, and~\ref{4root} that, in the case $K=\C$, this theorem means that the additive central charge~$c$ is always an integer, and an even integer if~$n$ is odd. The result is also consistent with a known fact about the Drinfel'd double: If $A=D(H)$ is the Drinfel'd double of an odd-dimensional semisimple Hopf algebra~$H$, then $n=\dim(H)^2$ is a square, and therefore congruent to~$1$ modulo~$4$. But in this case, it is known that \linebreak
$\chi_R(u) = \chi_R(u^{-1}) = \dim(H)$ (cf.~\cite{Mue2}, Prop.~5.18, p.~199; \cite{SoZhu}, Par.~6.1, p.~53).

However, we do not believe that the preceding theorem is the best possible result. Rather, we expect the following:
\begin{conj}
Consider a semisimple factorizable Hopf algebra~$A$ over a field~$K$ of characteristic zero. Let~$n$ be the dimension of~$A$, $\ua$ its Drinfel'd element, and~$\chi_R$ the character of its regular representation.
If~$n$ is even, then~$n$ is divisible by~$4$. Moreover, we have 
$\chi_R(u^{-1})^2 = \chi_R(u)^2$.
\end{conj}

Using the arguments of the proof of the preceding theorem together with Corollary~\ref{Cauchy}, the first part of this conjecture would follow from a conjecture of Y.~Kashina, who conjectured that the exponent~$N$ of~$A$ always divides the dimension~$n$ (cf.~\cite{KashAnti}, p.~1261). The best known result on this conjecture is a theorem by P.~Etingof and S.~Gelaki, which asserts that~$N$ divides~$n^3$ (cf.~\cite{EG4}, Thm.~4.3, p.~136). However, both Kashina's conjecture and our conjecture do not hold for quasi-Hopf algebras, as we will see in Paragraph~\ref{Semion} in an explicit example.

\subsection[Radford's example]{} \label{ExamRadf}
We have promised in Paragraph~\ref{ClassGauss} to explain why the Gaussian sum of a modular datum carries this name. Part of this explanation has
been given there already, when we saw that the Gaussian sum of a modular datum is often equal to the classical Gaussian sum. Here, we
now give an example of a modular datum whose Gaussian sum is exactly
the classical one; the name in the general case comes from this example.

Originally, this example was given by D.~E.~Radford (cf.~\cite{RadfAntiQuasi}, Sec.~3, p.~10; \cite{RadfKnotInv}, Sec.~2.1, p.~219); we have already contemplated it in this context in \cite{SoZhu}, Par.~5.5, p.~50, and we will use some of the computations carried out there. Consider a cyclic group~$G$ of order~$n$. Denote the group ring by~$A=K[G]$, and fix a generator~$g$ of~$G$. As~$A$ is cocommutative, $A$ is certainly quasitriangular with respect to the R-matrix~$\A \o \A$. However, with respect to this R-matrix, it is not factorizable. Radford has determined all possible R-matrices for~$A$, and shown that~$A$ can only be factorizable if~$n$ is odd, what we will assume for the rest of this paragraph, and that in this case the R-matrix necessarily has the form
$$R = \frac{1}{n} \sum_{i,j=0}^{n-1} \zeta^{-ij} g^i \o g^{j}$$
where~$\zeta$ is a primitive $n$-th root of unity (cf.~\cite{RadfKnotInv}, Sec.~2.3, p.~227).

The claim that the Gaussian sum~$\gp$ of the modular datum
that arises from this factorizable Hopf algebra is equal to the classical
Gaussian sum~$\Gp$ as introduced in Definition~\ref{ClassGauss} now follows directly from the formulas for the Drinfel'd element in this example (cf.~\cite{RadfKnotInv}, Sec.~2.1, p.~219; Sec.~2.3, p.~227; \cite{SoZhu}, Par.~5.5, p.~51): We get that
$$\chi_R(u^{-1}) = \Gp = \sum_{i=0}^{n-1} \zeta^{i^2}  \qquad \qquad
\chi_R(u) = \Gm = \sum_{i=0}^{n-1} \zeta^{-i^2}$$
Also, by using \cite{RadElemNum}, Chap.~11, Eq.~(11.14), p.~89, we observed in \cite{SoZhu}, Prop.~5.5, p.~51 that the Hopf symbol coincides with the Jacobi symbol in this example; a fact that motivated the definition of the Hopf symbol in the general case (cf.~\cite{SoZhu}, Def.~14.1, p.~114). As a consequence, we observed there
that by the first supplement to Jacobi's reciprocity law
we have in this example that
$$\chi_R(u^{-1}) =
\begin{cases}
\chi_R(u) &\text{if} \mspace{20mu} n \equiv 1 \pmod{4} \\
-\chi_R(u) &\text{if} \mspace{20mu} n \equiv 3 \pmod{4}
\end{cases}$$
Theorem~\ref{MainResult} should therefore be seen as a generalization of these two properties from this specific example to arbitrary odd-dimensional semisimple factorizable Hopf algebras.

\newpage
\section{Quasi-Hopf Algebras} \label{Sec:QuasiHopfAlg}
\subsection[Ribbon quasi-Hopf algebras]{} \label{RibQuasiHopf}
Theorem~\ref{MainResult} can be partially generalized to quasi-Hopf algebras. Although we will present the corresponding result, our main
focus is in a different direction: Quasi-Hopf algebras provide an example
of a normalized integral modular datum for which $\gp^4 = \gm^4$, but $\gp^2 \neq \gm^2$. In other words, this analogue of Theorem~\ref{MainResult}, and also Theorem~\ref{4root}, cannot be improved,
and the analogue of Conjecture~\ref{MainResult} for quasi-Hopf algebras is false.

Quasi-Hopf algebras were introduced by V.~G.~Drinfel'd in \cite{DrinfQuasiHopf}. We will use the setup from \cite{Kas}, Chap.~XV; however, we will assume that the unit constraints are trivial, i.e., that the elements denoted by~$r$ and~$l$ in \cite{Kas}, Prop.~XV.1.2, p.~369 are equal to~$1$. An algebra~$A$ over the field~$K$ is therefore called a quasi-Hopf algebra if it is equipped with algebra homomorphisms 
\mbox{$\da: A \rightarrow A^{\o 2} = A \o A$}, $\ea: A \rightarrow K$, and
$\sa: A \rightarrow A^{\op}$, called the coproduct, the counit, and the antipode, and elements $\Phi \in A^{\o 3}$, $\alpha \in A$, and 
$\beta \in A$ such that the axioms (1.1)--(1.4) in
\cite{Kas}, Prop.~XV.1.2, p.~369 as well as the axioms (5.1) and~(5.2) in
\cite{Kas}, Def.~XV.5.1, p.~379 are satisfied; in particular, $\Phi$ is assumed to be invertible. Note that these axioms imply
that~$A$ is an ordinary Hopf algebra if~$\Phi$, $\alpha$, and~$\beta$ are trivial in the sense that they are the unit elements in~$A^{\o 3}$ resp.~$A$. The element~$\Phi$ is called the associator of~$A$. We furthermore require that the antipode of~$A$ is bijective. We note that 
the antipode is compatible with the comultiplication in the sense that
$$\da(\sa(a)) = F^{-1} (\sa \o \sa)(\da^{\cop}(a)) F$$
(cf.~\cite{DrinfQuasiHopf}, Prop.~1.2, p.~1426), where~$F \in A^{\o 2}$
is an invertible element that is explicitly given in terms of the defining
structure elements (cf.~\cite{DrinfQuasiHopf}, Eq.~(1.36), p.~1429). Following the notational conventions of Paragraph~\ref{FactHopf}, we denote by~$F'$ the image of~$F$ under the interchange of the two tensor factors.

We also assume that~$A$ is quasitriangular. As in Paragraph~\ref{FactHopf},
this means that~$A$ is equipped with an
R-matrix~$R = \sum_i a_i \o b_i\in A^{\o 2}$;
however, the axioms given there must be modified as indicated in 
\cite{Kas}, Prop.~XV.2.2, p.~371. Note that the two sets of axioms coincide in the Hopf algebra case, where the associator is trivial. As explained in \cite{Kas}, p.~380, these structure elements can be used to turn the category of finite-dimensional left $A$-modules into an autonomous quasisymmetric category, where we, slightly deviating from \cite{Kas}, not require that autonomous categories be strict. 

The Drinfel'd element~$\ua$ that we considered in Paragraph~\ref{FactHopf}
can also be generalized to the quasi-Hopf algebra setting: It is then defined as
$$\ua = \sum_{i,j} \sa(b_i y_j \beta \sa(z_j)) \alpha a_i x_j$$
where we have used the notation 
$\Phi^{-1} = \sum_{j} x_j \o y_j \o z_j$. It has the following properties:
\begin{prop}
$\ua$ is invertible. Moreover, we have
\begin{enumerate}
\item
$\ea(\ua) = 1$

\item 
$\sa^2(a) = \ua a \ua^{-1}$

\item
$\da(\ua) = F^{-1} ((S \o S)(F')) (\ua \o \ua) (R'R)^{-1}$
\end{enumerate}
\end{prop}
\begin{pf}
The invertibility of~$\ua$ as well as the first and the second property are proved in \cite{AltCoste}, Sec.~3, p.~87f. The third property is also proved there (cf.~Eq.~(4.21), p.~95) under the assumption that~$\alpha$ is invertible. However, it was shown in \cite{BulNauRib}, Eq.~(3.6), p.~668 that this condition is unnecessary. Note that~$R'$ was defined in Paragraph~\ref{FactHopf}; the element~$F$ comes from the compatibility condition between coproduct and antipode mentioned above.
\qed
\end{pf}

The quasi-Hopf algebra~$A$ that we are considering will also be assumed to be a ribbon quasi-Hopf algebra, i.e., to be endowed with a ribbon element.
This means the following:
\begin{defn}
A nonzero central element $v \in A$ is called a ribbon element if it satisfies
$$\da(v) = (R'R) (v \o v) \quad \text{and} \quad \sa(v) = v$$
\end{defn}

We note that most authors use the inverse of this element as a ribbon element; for example, this is the convention in \cite{Kas}, Def.~XIV.6.1, p.~361. Our convention is the one used in \cite{Tur}, Sec.~XI.3.1, p.~500.
Ribbon elements are not unique: From a given ribbon element, we can construct another one by multiplying it with a central grouplike element that is invariant under the antipode. On the other hand,
this is obviously the only possible modification.

The basic properties of ribbon elements are given in the following corollary to the preceding proposition, which is proved in \cite{SoZhuRib}:
\begin{corollary}
A ribbon element~$v$ is invertible and satisfies~$\ea(v)=1$ as well as 
$$v^{-2} = \ua \sa(\ua)$$
\end{corollary}

This corollary in particular asserts that two of the axioms for ribbon quasi-Hopf algebras listed in \cite{BulNauRib}, Thm.~3.1, p.~667, namely Axiom~(3.1) and Axiom~(3.3), are actually consequences of the other axioms, namely Axiom~(3.2) and Axiom~(3.5). The axioms listed in \cite{BulNauRib} 
were an improved version of the original ones given in \cite{AltCoste},
Par.~4.1, p.~89. 

For an $A$-module~$V$, we define its twist~$\theta_V$ as
$$\theta_V: V \rightarrow V,~x \mapsto \theta_V(x) := vx$$
Because~$v$ is central, $\theta_V$ is $A$-linear, and therefore the
collection of all twists defines a natural transformation from the identity functor to itself. From this and the other axioms for a ribbon element, it follows that these twists turn the category of finite-dimensional $A$-modules into a ribbon category (cf.~\cite{Kas}, Def.~XIV.3.2, p.~349; \cite{Tur}, Sec.~I.1.4, p.~21). Accordingly, as already pointed out in Paragraph~\ref{ModCat}, every $A$-linear map $f: V \rightarrow V$ has a categorical trace~$\Tr_q(f)$. The categorical trace of~$f$ can be related to the ordinary trace by introducing the map
$$f_q: V \rightarrow V,~x \mapsto f(S(\alpha) u v \beta x)$$
As explained in \cite{AltCoste}, Par.~4.4, p.~95, the ordinary trace of this map is equal to the categorical trace of~$f$; in other words:
\begin{lemma}
$\Tr_q(f) = \Tr(f_q)$
\end{lemma}
We also note that the version of this lemma for ordinary Hopf algebras, which we have already used in Paragraph~\ref{MainResult}, 
can be found in \cite{Kas}, Prop.~XIV.6.4, p.~363 and
\cite{Tur}, Lem.~XI.3.3, p.~501.

\subsection[Involutory quasi-Hopf algebras]{} \label{InvolQuasiHopf}
On our ribbon quasi-Hopf algebra~$A$, we now impose a number of additional 
assumptions: First, we assume that the base field~$K$ is algebraically closed of characteristic zero. Second, we assume that~$A$ is of finite dimension~$n$ and semisimple as an algebra. As for ordinary Hopf algebras
in Paragraph~\ref{FactHopf}, the Wedderburn structure theorem then provides 
us with a family of simple modules $(V_i)_{i \in I}$, one of which, denoted~$V_o$, is equal to the base field~$K$, turned into an $A$-module via the counit. The category of finite-dimensional $A$-modules is therefore a semisimple ribbon category; it will be modular if the Verlinde matrix~$\V$
from Paragraph~\ref{ModCat} is invertible, as we from now on assume.
Finally, we assume that $v^{-1} = S(\alpha) u \beta$. We will comment below on the meaning of this assumption; for the moment, we only note that it implies that the categorical trace agrees with the
ordinary trace, because we then have~$f=f_q$ in Lemma~\ref{RibQuasiHopf}. As in Paragraph~\ref{FactHopf}, the entries of the Verlinde matrix~$\V$ are therefore given as
$$s_{ij} = (\chi_i \o \chi_{j^*})(R'R)$$
where $\chi_i$ denotes the character of~$V_i$ and~$i^* \in I$ is the unique index satisfying $V_{i^*} \cong V_i^*$.

The remaining parts of the discussion in Paragraph~\ref{FactHopf} also 
carry over: The diagonal entry~$t_i$
of the Dehn matrix~$\T$ is determined by the condition that~$v$
acts on~$V_i$ as $t_i \id_{V_i}$, so that we arrive again at a normalized modular datum $(I,o,{}^*,\V,\T)$. Using the same computation as in Paragraph~\ref{FactHopf}, we get that $n_i = s_{io} = \dim(V_i)$, the dimension of a vector space, so that this modular datum is again integral.
As before, the Wedderburn structure theorem implies that the global dimension~$n$ of the modular datum is equal to the dimension of~$A$, and the exponent~$N$ of the modular datum is again equal to the order of the ribbon element~$v$. Another consequence of the Wedderburn structure theorem is
that the character of the regular representation is given by the formula~$\chi_R(a) = \sum_{i \in I} n_i \chi_i(a)$, so that the Gaussian sum $\gp = \sum_{i \in I} n_i^2 t_i$ becomes $\gp = \chi_R(v)$,
and the reciprocal Gaussian sum becomes $\gm = \chi_R(v^{-1})$,
exactly as in Paragraph~\ref{FactHopf}, where we had $v=u^{-1}$.
We can therefore express the fusion symbol~$\fs(q)$ again in terms of powers of the ribbon element, and also reformulate Proposition~\ref{DefModData}.\ref{DefModDataRecGauss} as the statement $\chi_R(v) \chi_R(v^{-1})= \dim(A)$.

Let us now try to put our assumption $v^{-1} = S(\alpha) u \beta$
into context. We first note that this assumption has the following consequence:
\begin{prop}
$A$ is involutory.
\end{prop}
\begin{pf}
By definition (cf.~\cite{BulCaenTorInvolut}, Def.~3.1, p.~263), we have to
show that
$$\sa^2(a) = (\sa(\beta) \alpha) \; a \; (\beta \sa(\alpha))$$
By assumption, we have
$v^{-1} = \sa(\alpha) \ua \beta = \sa(\alpha) \sa^2(\beta) \ua$
and therefore 
$$\sa(\beta) \alpha = \sa^{-1}(v^{-1}\ua^{-1}) = \sa^{-1}(\ua^{-1}) v^{-1}
= \sa(\ua^{-1}) v^{-1}$$
Since~$A$ is finite-dimensional, this shows that~$\alpha$ and~$\beta$ are invertible. From Proposition~\ref{RibQuasiHopf}, we also have 
$\sa^2(a) = \sa(\ua)^{-1} a  \sa(\ua)$. Because~$v$ is central, this implies
\begin{align*}
(\sa(\beta) \alpha) \; a = \sa(\ua^{-1}) v^{-1} \;  a 
= \sa^2(a) \;  \sa(\ua^{-1}) v^{-1} = \sa^2(a) \; (\sa(\beta) \alpha)
\end{align*}
Therefore, our assertion will follow if we can show that
$\sa(\beta) \alpha = (\beta \sa(\alpha))^{-1}$. Solving the assumption
$v^{-1} = \sa(\alpha) \ua \beta$ for~$\ua$, we get 
$$\ua = \sa(\alpha)^{-1} v^{-1} \beta^{-1} = 
\sa(\alpha)^{-1} \beta^{-1} v^{-1}$$
But by Corollary~\ref{RibQuasiHopf}, we have $v^{-2} = \ua \sa(\ua)$,
so that
$$\sa(\alpha)^{-1} \beta^{-1} = \ua v = \sa(\ua^{-1}) v^{-1}
= \sa(\beta) \alpha$$
by the first computation above.
\qed
\end{pf}

The nature of our assumption becomes clearer when restated in categorical terms: In every left rigid quasisymmetric monoidal category, i.e., a 
quasisymmetric monoidal category with left duality, there is 
a canonical natural transformation between the identity functor
and the double dual functor~$V \mapsto V^{**}$. In the case of 
the category of finite-dimensional modules over a quasitriangular
quasi-Hopf algebra~$A$, this is the natural transformation
$V \rightarrow V^{**},~x \mapsto \Theta(ux)$, where~$\Theta$
is the isomorphism between~$V$ and its double dual from linear algebra,
defined as $\Theta(x)(\varphi) = \varphi(x)$ for $x \in V$ and
$\varphi \in V^*$. However, this is not a monoidal transformation,
but rather a ribbon transformation (cf.~\cite{So5}, Def.~3.2, p.~439),
which is reflected by the identities in Proposition~\ref{RibQuasiHopf}.
Monoidal transformations between the identity functor
and the double dual functor can also be obtained from elements~$w \in A$;
however, in comparison to Proposition~\ref{RibQuasiHopf}, such an element $w$ should satisfy
\begin{enumerate}
\item
$\ea(w) = 1$

\item 
$\sa^2(a) = w a w^{-1}$

\item
$\da(w) = F^{-1} ((S \o S)(F')) (w \o w)$
\end{enumerate}
(cf.~\cite{BulCaenTorInvolut}, Prop.~4.2, p.~273; \cite{MasNg}, Thm.~7.1, p.~182; \cite{SchauenQuasiHopf}, Rem.~3.3, p.~133). We then get a monoidal transformations between the identity functor and the double dual functor
by assigning to every object~$V$ the morphism $V \rightarrow V^{**},~x \mapsto \Theta(wx)$. If~$w$ is invertible, then this will be a natural equivalence, i.e., a pivotal structure. It is now obvious from Proposition~\ref{RibQuasiHopf} and the definition of a ribbon element that the element~$w:=uv$ satisfies these equations. Taking a second look at the preceding proof, we see that our condition $v^{-1} = S(\alpha) u \beta$ is equivalent to the requirement that $w = \sa(\beta)\alpha$; in other words, we are using the unique pivotal structure for which categorical traces and ordinary traces coincide (cf.~\cite{BulCaenTorInvolut}, Prop.~4.3, p.~275; \cite{EtNikOst}, Prop.~8.23, Prop.~8.24, p.~622f).

\subsection[Results for quasi-Hopf algebras]{} \label{ResultQuasiHopf}
For a quasi-Hopf algebra~$A$ that satisfies the assumptions in Paragraph~\ref{InvolQuasiHopf}, the methods developed in Section~\ref{Sec:OddExp} now apply in the same way as they applied to an ordinary Hopf algebra in Paragraph~\ref{MainResult}, and we get the following result:
\begin{thm}
If the dimension~$n$ of~$A$ is odd, we have
$$\chi_R(v^{-1}) = 
\begin{cases}
\chi_R(v) &: n \equiv 1 \pmod{4} \\
- \chi_R(v) &: n \equiv 3 \pmod{4}
\end{cases}$$
Moreover, the fusion symbol coincides with the Jacobi symbol; i.e.,
we have $\fs(q) = \jac{q}{n}$ for all~$q \in \Z$.
\end{thm}
\begin{pf}
By Cauchy's theorem for quasi-Hopf algebras (cf.~\cite{NgSchauen3}, Thm.~8.4, p.~63), $N$ and~$n$ have the same prime divisors, so~$N$ is also odd. Note that it follows from~\cite{NgSchauen3}, Thm.~7.7, p.~60 that~$N$
is the Frobenius-Schur exponent of~$A$. Therefore, the equation $\chi_R(v^{-1}) = (-1)^{(n-1)/2} \chi_R(v)$ follows again directly from Theorem~\ref{ClassGauss}. Cauchy's theorem also implies that an integer~$q$ is relatively prime to~$N$ if and only if it is relatively prime to~$n$, which means that $\fs(q) \neq 0$ if and only if $\jac{q}{n} \neq 0$. But as soon as these two numbers are nonzero, they are equal by Theorem~\ref{ClassGauss}.
\qed
\end{pf}

It is difficult to overlook that the preceding theorem does not contain
an analogue of the second part of Theorem~\ref{MainResult}, which was
concerned with the case where the dimension~$n$ is even, despite 
the fact that the projective congruence subgroup theorem has already been
carried over from Hopf algebras (cf.~\cite{SoZhu}, Thm.~9.4, p.~94) to quasi-Hopf algebras (cf.~\cite{NgSchauen4}, Thm.~8.8, p.~35). The reason for this is that the Galois property (cf.~\cite{SoZhu}, Lem.~12.2, p.~115)
has not yet been carried over. However, we expect that this property can
also be established for quasi-Hopf algebras in an analogous way, so that then the methods of Section~\ref{Sec:EvenExp} would apply and we would
get that $\chi_R(v^{-1})^4 = \chi_R(v)^4$ for any quasi-Hopf algebra that 
satisfies the conditions in Paragraph~\ref{InvolQuasiHopf}.

On the other hand, we have explicitly formulated in Conjecture~\ref{MainResult} our expectation that, in the Hopf algebra
case, we even have the stronger result that 
$\chi_R(v^{-1})^2 = \chi_R(v)^2$. This result, if correct, at least does not
carry over to quasi-Hopf algebras, as we will show now by an explicit example.

\subsection[3-cocycles]{} \label{3Cocyc}
In general, if~$G$ is a finite group and 
$\omega: G \times G \times G \rightarrow K^\times$ is a normalized 3-cocycle (cf.~\cite{Kas}, Sec.~XV.5, Eq.~(5.3), p.~380), then the dual 
group ring~$K[G]^*$, which is an ordinary Hopf algebra, can also be considered as a quasi-Hopf algebra with respect to the
associator
$$\Phi = \sum_{g,h,k \in G} \omega(g,h,k) \; e_g \o e_h \o e_k$$ 
and the antipode elements
$$\alpha = 1 \qquad  \qquad 
\beta = \sum_{g \in G} \omega(g^{-1},g,g^{-1}) \; e_g$$
where $e_g$ is the dual basis element of~$g \in G \subset K[G]$
(cf.~\cite{HausNillDoub}, App.~A, p.~585; \cite{NgSchauen2}, Sec.~7, p.~1856). In particular, if $G=\Z_n$ is the cyclic group of order~$n$, which we represent in the form $\Z_n:=\{0,1,2,\ldots,n-1\}$, we can
construct such a cocycle in the following way: 
For $i,j \in \{0,1,2,\ldots,n-1\}$, we define
$$q_{ij} := \frac{1}{n} (\bar{i} + \bar{j} - \overline{i+j})
= \begin{cases}
0 &: i+j \le n-1\\ 1 &: i+j \ge n
\end{cases}$$
where, for $i \in \Z$, the element $\bar{i} \in \{0,1,2,\ldots,n-1\}$ is the unique number that satisfies $i \equiv \bar{i} \pmod{n}$. 
Then the function
$$\sigma: \Z_n \times \Z_n \rightarrow K^\times,~(i,j) \mapsto
\sigma(i,j) := \zeta^{q_{ij}}$$
where~$\zeta$ is an $n$-th root of unity, is a normalized 2-cocycle
(cf.~\cite{Jac3}, Chap.~I, \S~15, p.~80), and from this it is comparatively easy to see that
$$\omega: \Z_n \times \Z_n \times \Z_n \rightarrow K^\times,~(i,j,k) \mapsto \omega(i,j,k) := \sigma(i,j)^k$$ is a normalized 3-cocycle.
Further details on this discussion and related material can be found in
\cite{MoSeib}, Eq.~(E.14), p.~251;
\cite{NatExpFin}, Ex.~5.8, p.~265; and \cite{NgSchauen3}, Sec.~9, p.~68f.

\subsection[The two-dimensional case]{} \label{TwoDim}
The simplest nontrivial case of the above construction is the case~$n=2$.
This special case has already been used as an example several times, as in \cite{EGQuasiHopf}, Par.~2.3, p.~687 and \cite{NgSchauen2}, Ex.~5.4, p.~1854. Note that the convention for~$\alpha$ and~$\beta$ in \cite{EGQuasiHopf} slightly deviates from the one used here. As we want $\Phi \neq 1$, we have to
choose $\zeta = -1$; we then have
$$q_{ij} = 
\begin{cases}
0 &: i=0 \; \text{or} \; j=0 \\
1 &: i=1 \; \text{and} \; j=1
\end{cases}
$$
so that $q_{ij} = ij$ and consequently
$$\sigma(i,j) = (-1)^{ij} \qquad \qquad \omega(i,j,k) = (-1)^{ijk}$$
The possible R-matrices for this quasi-Hopf algebra have been determined
in \cite{BulCaenTorInvolut}, Prop.~3.6, p.~265; they have the form
$$R = e_0 \o e_0 + e_0 \o e_1 + e_1 \o e_0 + \iota e_1 \o e_1$$
where $\iota \in K$ is a primitive fourth root of unity. Inserting this into the definition in Paragraph~\ref{RibQuasiHopf}, we find that the  Drinfel'd element then is $\ua = e_0 + \iota e_1$. We also compute the element~$F$ that we introduced in Paragraph~\ref{RibQuasiHopf}:
\begin{lemma}
$$F = e_0 \o e_0 + e_0 \o e_1 + e_1 \o e_0 - e_1 \o e_1$$
\end{lemma}
\begin{pf}
We have $\alpha = 1 = e_0 + e_1$ and $\beta = e_0 - e_1$.
It therefore follows from \cite{DrinfQuasiHopf}, Eq.~(1.35), p.~1429
that~$F$ is in our case equal to the element~$\gamma$ defined in  \cite{DrinfQuasiHopf}, Eq.~(1.24), p.~1426.
Now we have
\begin{align*}
&(1 \o \Phi^{-1}) (\id \o \id \o \id \o \Delta)(\Phi) = \\
&(\sum_{p,q,r =0}^1 (-1)^{pqr} \; 1 \o e_p \o e_q \o e_r)
(\sum_{i,j,k,l =0}^1 (-1)^{ij(k+l)} \; e_i \o e_j \o e_k \o e_l)= \\
&\sum_{i,j,k,l =0}^1 
(-1)^{jkl} (-1)^{ij(k+l)} \; e_i \o e_j \o e_k \o e_l
\end{align*}
Inserting this into the definition of~$\gamma$, we find
\begin{align*}
\gamma &= \sum_{i,j,k,l =0}^1 
(-1)^{jkl} (-1)^{ij(k+l)} \; S(e_j) \alpha e_k \o S(e_i) \alpha e_l\\
&= \sum_{i,j =0}^1 
(-1)^{j^2i} (-1)^{ij(j+i)} \; e_j \o e_i 
= \sum_{i,j =0}^1 (-1)^{ij} \; e_j \o e_i
\end{align*}
as asserted.
\qed
\end{pf}

As a consequence, Proposition~\ref{RibQuasiHopf}.3 implies 
that~$\ua^{-1}= e_0 - \iota e_1$ is a ribbon element. However, it is not the only one: As we pointed out directly after Definition~\ref{RibQuasiHopf}, we can modify a given
ribbon element by multiplying it with a central grouplike element that
is invariant under the antipode. Therefore,
$v := \beta \ua^{-1} = e_0 + \iota e_1 = \ua$ is another ribbon element,
and this is the one we want, because we want to satisfy the requirement
$v^{-1} = \sa(\alpha) \ua \beta$ made in Paragraph~\ref{InvolQuasiHopf}.

Let us consider how we can satisfy the other requirements made in Paragraph~\ref{InvolQuasiHopf}. Using $I = \{0,1\}$ as index set
with $o=0$ as distinguished element, we have two simple modules
$V_0 := K \cong K e_0$ and $V_1 := K e_1$. Note that $K e_0$ and
$K e_1$ are ideals. The characters~$\chi_0$ and~$\chi_1$ of these modules 
satisfy $\chi_i(e_j)=\delta_{ij}$. Because the antipode is the identity in this example, the involution~${}^*$ is also the identity here.
The formula $s_{ij} = (\chi_i \o \chi_{j^*})(R'R)$ for the entries of the 
Verlinde matrix given in Paragraph~\ref{InvolQuasiHopf} therefore implies that this matrix is
$$\V = \begin{pmatrix} 1 & 1 \\ 1 & -1\end{pmatrix}$$
so that the last requirement, the invertibility of~$\V$, is also met.
The form of the Dehn matrix follows directly from the form of the ribbon element given above; we have
$$\T = \begin{pmatrix} 1 & 0 \\ 0 & \iota \end{pmatrix} $$

\subsection[The semion datum]{} \label{Semion}
The modular datum that we have just obtained from the quasi-Hopf algebra
$K[\Z_2]^*$ also arises from other constructions, for example from the Kac-Moody algebra~$A_1^{(1)}$ at level~$1$ (cf.~\cite{GanModDat}, Sec.~3, Eq.~(3.5), p.~223; \cite{GanMoonMonst}, Sec.~6.2.1, Eq.~(6.2.2), p.~368; \cite{K2}, Chap.~13, Ex.~13.8, p.~265). Following \cite{RowStongWang}, Sec.~5.3, p.~38, we call this datum the semion datum. It is frequently used in the literature; we mention only \cite{EtVafa}, Sec.~5, p.~655 and \cite{Hu}, Par.~4.1.3, p.~56.

The semion datum is a normalized integral modular datum, like every modular datum that arises from a quasi-Hopf algebra that satisfies the requirements
imposed in Paragraph~\ref{InvolQuasiHopf}. Of course, it is also easy
to verify the axioms in Definition~\ref{DefModData} directly. Obviously,
the exponent, which is equal to the normalized exponent here, is~$N=N_o=4$,
and the global dimension~$n$ is~$2$, the dimension of the quasi-Hopf algebra. For the Gaussian sum, we find~$\gp = 1+\iota$, and for the reciprocal Gaussian sum~$\gm = 1-\iota$, so that we have
$\gp^2 = 2\iota = - \gm^2$. This shows that Theorem~\ref{4root}
cannot be improved, because the semion datum satisfies the hypotheses
of that theorem:
\begin{lemma}
The semion datum is a projective congruence datum that is also Galois.
\end{lemma}
\begin{pf}
The Galois group of the cyclotomic field~$\Q_4 \subset K$ contains only
one nonidentity element, namely the automorphism~$\gamma$ discussed in Paragraph~\ref{ComplConj}. The Galois condition therefore becomes
$t_{\gamma.i} = \gamma^2(t_i)$, where $t_0=1$ and~$t_1=\iota$. But by
Proposition~\ref{ComplConj}, we have $\gamma.i = i^* =i$, so that the condition holds.

As we already mentioned in Paragraph~\ref{ResultQuasiHopf}, the projective congruence subgroup theorem for quasi-Hopf algebras implies that the
semion datum is a projective congruence datum. In the present case,
however, this is easy to see directly, because to obtain~$\SL(2,\Z_4)$ from~$\SL(2,\Z)$, the only defining
relation that we need to add is the relation~$\bar{\gt}^4=1$ (cf.~\cite{Fine}, Sec.~3.6, p.~64; \cite{Newman}, Chap.~VIII, p.~145ff), which is clearly satisfied in our example.
\qed
\end{pf}

However, the projective representation
$$\SL(2,\Z_4) \rightarrow \PGL(2,K),~\bar{\gv} 
\mapsto \V,~\bar{\gt} \mapsto \T$$
that we have just described cannot be lifted to an ordinary representation
of~$\SL(2,\Z_4)$. To see this, assume that there is such a representation,
and denote the images of the generators~$\bar{\gv}$ and~$\bar{\gt}$ by~$\V'$ and~$\T'$.
These matrices then have to satisfy the relations
$$\V'^4 = \E \qquad \qquad \T'^4 = \E \qquad \qquad (\T' \V')^3 = \V'^2$$
Because the representation lifts the projective one, there must be nonzero
numbers $D,\ell \in K^\times$ such that
$$\V' = \frac{\V}{D} \qquad \text{and} \qquad \T' = \frac{\T}{\ell}$$
Inserting this into the preceding relations, these become
$$\V^4 = D^4 \E \qquad \qquad \T^4 = \ell^4 \E \qquad \qquad 
(\T \V)^3 = D \ell^3 \V^2$$
But we know already that $\T^4 = \E$, and from Axiom~\ref{DefModData}.\ref{DefModData3} and the constant form of
Axiom~\ref{DefModData}.\ref{DefModData4} we know that
$\V^2 = 2 \E$ and $(\T \V)^3 = \gp \V^2$. Therefore, our conditions imply
that 
$$D^4 = 4 \qquad \qquad \ell^4 = 1 \qquad \qquad 
D \ell^3 = \gp = 1 + \iota$$
Raising the last condition to the fourth power, we get
$$4 = (D \ell^3)^4 = (1 + \iota)^4 = (2\iota)^2 = -4$$
a contradiction.

Stated differently, this computation shows that it is not possible to
extend the semion datum by choosing a generalized rank~$D$ and a (multiplicative) central charge~$\ell$ in such a way that the extended
modular datum is a congruence datum. However, the kernel of the linear representation of the modular group that we have associated with any extension in Paragraph~\ref{CentCharge} is still
a congruence subgroup, although not of level~$N=4$:
\begin{prop}
It is possible to extend the semion datum in such a way that the kernel of the associated linear representation of the modular group contains~$\Gamma(8)$. For every extension, this kernel contains~$\Gamma(24)$.
\end{prop}
\begin{pf}
\begin{list}{(\arabic{num})}{\usecounter{num} \leftmargin0cm \itemindent5pt}
\item
Using the abbreviation 
$\bar{\gd} = \bar{\gd}(5,5)= 
\bar{\gv} \bar{\gt}^{5} \bar{\gv}^{-1} \bar{\gt}^5 \bar{\gv} \bar{\gt}^{5}$ from Definition~\ref{DiagMat}, the matrices 
$\bar{\gv}, \bar{\gt} \in \SL(2,\Z_8)$ satisfy
$$\bar{\gt}^8 = 1 \qquad \qquad  \bar{\gv} = \bar{\gd} \bar{\gv} \bar{\gd} \qquad \qquad  
\bar{\gt} = \bar{\gd} \bar{\gt} \bar{\gd}^{-1}$$
Together with the relations $\bar{\gv}^4 = 1$ and 
$(\bar{\gt} \bar{\gv})^3 = \bar{\gv}^2$ from Paragraph~\ref{CentCharge},
these are defining relations for~$\SL(2,\Z_8)$. This follows by transposing the presentation given in \cite{Men}, Lem.~3.3, p.~210 and noting that the
last relation there is equivalent to Equation~(3.17) on the following page there.

Now let~$D$ be one of the two solutions to the equation~$D^2 = -2$,
and set $\ell:=(1-\iota)/D$. Then~$D$ is a generalized rank that is not a rank. Furthermore, we have $\ell^2 = -2\iota/D^2 = \iota$ and therefore
$\ell^3 = (\iota+1)/D = \gp/D$. This shows that~$\ell$ is a multiplicative
central charge, so that we get an extended modular datum. As we have discussed in Paragraph~\ref{CentCharge}, the associated homogeneous Verlinde matrix~$\V'=\V/D$ and the homogeneous Dehn matrix~$\T'=\T/\ell$
give rise to a linear representation of the modular group by mapping~$\gv$ to~$\V'$ and~$ \gt$ to~$\T'$. The kernel of this representation will contain~$\Gamma(8)$ if the additional relations stated above are satisfied.
Because $\ell^2 = \iota$, we have $\ell^8 = 1$, so $\T'^8 = \E$, which means that the first relation holds.

Because $\iota^5 = \iota$, we have $\T^5= \T$, and therefore the image
of~$\gd(5,5)$ under the associated linear representation of the modular group is
$$\frac{1}{D \ell^{15}} \V \T \V^{-1} \T \V \T = 
\frac{1}{2 D \ell^{15}} (\V \T)^3 
= \frac{\gp}{2 D \ell^{15}} \V^2
= \frac{\gp}{D \ell^{15}} \E = \frac{\gp \ell}{D} \E = \ell^4 \E = - \E$$
and from this it is obvious that the two remaining relations are also satisfied.

\item
To prove that, for any extension, $\Gamma(24)$ is in the kernel of the associated representation of the modular group, recall that we have 
$\Z_{24} \cong \Z_8 \times \Z_3$ by the Chinese remainder theorem, and therefore $\SL(2,\Z_{24}) \cong \SL(2,\Z_8) \times \SL(2,\Z_3)$.
The group $\SL(2,\Z_8) \times \SL(2,\Z_3)$
is easy to present: It is generated by the four elements
$(\bar{\gt},1)$, $(\bar{\gv},1)$, $(1,\bar{\gt})$, and $(1,\bar{\gv})$.
The generators $(\bar{\gt},1)$ and $(\bar{\gv},1)$ have to satisfy the 
preceding relations for $\SL(2,\Z_8)$. The generators $(1,\bar{\gt})$ and $(1,\bar{\gv})$ have to satisfy the defining relations for $\SL(2,\Z_3)$, which are those for the modular group stated in Paragraph~\ref{CentCharge}
and the relation $(1,\bar{\gt})^3=(1,1)$ (cf.~\cite{Fine}, Sec.~3.6, p.~64; \cite{Newman}, Chap.~VIII, p.~145ff). Finally, each of the first two generators has to commute with each of the latter two.

This presentation translates into a presentation of $\SL(2,\Z_{24})$
if we can find the corresponding generators. For this, note that the isomorphism coming from the Chinese remainder theorem maps
$\bar{9} \in \Z_{24}$ to~$(\bar{1},\bar{0}) \in Z_8 \times \Z_3$,
and similarly~$\bar{16}$ to~$(\bar{0},\bar{1})$. As a consequence,
$\bar{\gt}^9 \in \SL(2,\Z_{24})$ is mapped to the first generator
$(\bar{\gt},1) \in \SL(2,\Z_8) \times \SL(2,\Z_3)$, and similarly
$\bar{t}^{16}$ is mapped to~$(1,\bar{\gt})$.

To find the preimages of the two remaining generators, note that it follows
from Definition~\ref{DiagMat} that $\gd(0,0) = \gv$ and $\gd(1,1) = 1$.
Consequently, the element $\bar{\gd}(16,16) \in \SL(2,\Z_{24})$ is mapped to
$(\bar{\gv},1) \in \SL(2,\Z_8) \times \SL(2,\Z_3)$, and similarly
$\bar{\gd}(9,9)$ is mapped to $(1,\bar{\gv})$.

\item
Now let~$D$ be an arbitrary generalized rank, and~$\ell$ be an arbitrary multiplicative central charge. According to Definition~\ref{CentCharge},
this means that $D^4=n^2=4$ and $\ell^3 = \gp/D$, which implies that
$$\ell^{12} = \frac{\gp^4}{D^4} = \frac{1}{4} (1+\iota)^4 = -1$$
As explained in Paragraph~\ref{CentCharge}, the associated linear representation of the modular group maps~$\gv$ to the homogeneous
Verlinde matrix~$\V'=\V/D$ and~$\gt$ to the homogeneous Dehn matrix~$\T'=\T/\ell$. If this representation would factor over $\SL(2,\Z_{24})$, it would therefore map the generator $\bar{\gt}^9$
to~$\T'^9=\T/\ell^9$ and the generator $\bar{\gt}^{16}$
to $\T'^{16}=\E/\ell^{16}$. Using this, we see that it would map the generator $\bar{\gd}(9,9)$ to~$\E$ and the generator
$\bar{\gd}(16,16)$ to~$\V'=\V/D$.

This factorization will occur if these potential images satisfy the relations that we have described above. Since the generators~$\bar{\gt}^{16}$
and~$\bar{\gd}(9,9)$, which correspond to~$(1,\bar{\gt})$ and~$(1,\bar{\gv})$, are mapped to multiples of the unit matrix, the commutation property described above is satisfied. The only remaining conditions that these two generators 
have to satisfy are the defining relations of~$\SL(2,\Z_3)$ already outlined above, which here reduce to the condition~$(\ell^{16})^3=1$,
an obvious consequence of~$\ell^{12}=-1$.

\item \enlargethispage{-1mm}
It remains to check that the images~$\V/D=\V'$ and~$\T/\ell^9=\T'/\ell^8$ of the generators that correspond to~$(\bar{\gv},1)$ and~$(\bar{\gt},1)$ satisfy the defining relations of $\SL(2,\Z_8)$ described in the first step.
Because the relations in~$\SL(2,\Z)$ hold, we know that
$\V'^4=\E$, and $(\T'\V')^3=\V'^2$, and these relations still hold if
we rescale~$\T'$ by~$1/\ell^8$ and leave~$\V'$ untouched. But as we saw in the first step, there are three additional conditions. The first one
is that the order of~$\T/\ell^9$ divides~$8$, which holds because
$\T^4=1$ and $(\ell^9)^8=(\ell^{12})^6=(-1)^6=1$. For the remaining
two conditions, we need to find, as in the first step, the potential image of~$(\bar{\gd}(5,5),1)$. This is found by a very similar computation, 
namely
$$\frac{1}{D \ell^{135}} \V \T \V^{-1} \T \V \T = 
\frac{1}{2 D \ell^{135}} (\V \T)^3 
= \frac{\gp}{2 D \ell^{135}} \V^2
= \frac{\gp}{D \ell^{135}} \E = - \E$$
where we have used that 
$((\ell^9)^5)^3 = \ell^{135} = (\ell^{12})^{11} \ell^3 =
(-1)^{11} \ell^3 = - \ell^3 = - \gp/D$.
As in the first step, the two remaining conditions are immediate consequences of this fact.
\qed
\end{list}
\end{pf}

It is instructive to compare our treatment of the semion datum with the one in~\cite{Eholz1} and~\cite{RowStongWang}. For this, we again 
consider the base field~$K=\C$ and write, as in Paragraph~\ref{CentCharge}, $\ell=e^{2\pi ic/24}$, where~$c$ is the additive central charge. 
Working with the homogeneous versions, Eholzer first considers in~\cite{Eholz1}, Sec.~5.1, p.~638 the matrices
$$\V' = \frac{1}{\sqrt{2}}\begin{pmatrix} -1 & -1 \\ -1 & 1 \end{pmatrix}
\qquad \qquad 
\T' = \begin{pmatrix} e^{2\pi i/8} & 0 \\ 0 & e^{2\pi i(3/8)} \end{pmatrix}
= \frac{1}{e^{-2 \pi i/8}}\begin{pmatrix} 1 & 0 \\ 0 & e^{\pi i/2} \end{pmatrix}$$
Comparing with our language, we see that Eholzer uses~$D=-\sqrt{2}$, which is a rank. Furthermore, we see that we have $\iota = e^{\pi i/2} = i$, and 
$\ell = e^{-2 \pi i/8} = e^{2 \pi i(21/24)}$, so  that the additive central charge is~$c=21$. Although these choices for~$\ell$ and~$D$ differ from the ones that we made in the first step of the preceding proof, these choices also entail that $D\ell^3=1+\iota$, and $\ell^4=-1$, and a second look at the preceding argument shows that this implies that~$\Gamma(8)$ lies in the kernel, which is in fact the starting point for Eholzer's construction of the example. The second version of the homogeneous Dehn matrix that Eholzer considers is the complex conjugate of the first, for which we therefore have~$c=3$.

On the other hand, Rowell et al.\ consider in~\cite{RowStongWang}, Sec.~5.3, p.~38 the matrices
$$\V = \begin{pmatrix} 1 & 1 \\ 1 & -1 \end{pmatrix}
\qquad \qquad 
\T = \begin{pmatrix} 1 & 0 \\ 0 & i \end{pmatrix}$$
Note that our matrix~$\V$ corresponds to their matrix~$\tilde{S}$ and
our matrix~$\V'$ corresponds to their matrix~$S$. They then use
$D=\sqrt{2}$, the other rank, and therefore get from their analogue of Definition~\ref{CentCharge} (cf.~\cite{RowStongWang}, Prop.~2.6, p.~9)
directly~$c=1$ for the central charge. In our terminology, this leads to
the choice~$\ell=e^{2 \pi i/24}$, so that we do not have~$\ell^8=1$. Accordingly, the kernel of the associated linear representation of the modular group does not contain~$\Gamma(8)$, but only~$\Gamma(24)$. This 
distinction is, however, not relevant for the authors, because they only
look at the projective representation (cf.~\cite{RowStongWang}, p.~44 and p.~47).

\newpage
\addcontentsline{toc}{section}{Bibliography} \label{Bibliography}


\begin{thebibliography}{99}
\bibitem{AltCoste} D.~Altsch\"uler/A.~Coste: Quasi-quantum groups, knots, three-manifolds, and topological field theory, Commun.~Math.~Phys.~150 (1992), 83-107

\bibitem{Apos} T.~M.~Apostol: Modular forms and Dirichlet series in number theory, 2nd ed., Grad.~Texts Math., Vol.~41, Springer, Berlin, 1990

\bibitem{BakKir} B.~Bakalov/A.~Kirillov Jr.: Lectures on tensor categories and modular functors, Univ.~Lect.~Ser., Vol.~21, Am.~Math.~Soc., Providence, 2001

\bibitem{BanModGal} P.~Bantay: The kernel of the modular representation and the Galois action in RCFT, Commun.\ Math.~Phys.~233 (2003), 423-438

\bibitem{Beauv} A.~Beauville: Conformal blocks, fusion rules and the Verlinde formula. In: M.~Teicher (ed.): Proceedings of the Hirzebruch 65 conference on algebraic geometry, Isr.~Math.~Conf.~Proc., Vol.~9, Bar-Ilan Univ., Ramat-Gan, 1996,~75-96

\bibitem{BeEvWil} B.~C.~Berndt/R.~J.~Evans/K.~S.~Williams: Gauss and Jacobi sums, Can.\ Math.\ Soc.\ Ser.\ Monogr.\ Adv.\ Texts, Vol.~21, Wiley, New York, 1998

\bibitem{Beyl} F.~R.~Beyl: The Schur multiplicator of~$\SL(2,\Z/m\Z)$
and the congruence subgroup property, Math.~Z.~191 (1986), 23-42

\bibitem{BulCaenTorInvolut} D.~Bulacu/S.~Caenepeel/B.~Torrecillas:
Involutory quasi-Hopf algebras, Algebr.\ Represent.\ Theory~12 (2009), 257-285

\bibitem{BulNauRib} D.~Bulacu/E.~Nauwelaerts: Quasitriangular and ribbon
quasi-Hopf algebras, Commun.~Algebra~31 (2003), 657-672

\bibitem{ChariPress} V.~Chari/A.~Pressley: A guide to quantum groups, Camb.~Univ.~Press, Cambridge, 1994

\bibitem{CostGann} A.~Coste/T.~Gannon: Congruence subgroups and rational conformal field theory, Preprint, math.QA/9909080, 1999

\bibitem{CostGannRue} A.~Coste/T.~Gannon/P.~Ruelle: Finite group modular data, Nuclear Phys.~B~581 (2000), 679-717

\bibitem{CoxMos} H.~S.~M.~Coxeter/W.~O.~J.~Moser: Generators and relations for discrete groups, 4th~ed., Ergeb.\ Math.\ Grenzgeb., Vol.~14, Springer, Berlin, 1984

\bibitem{Cuntz} M.~Cuntz: Integral modular data and congruences, J.~Algebr.~Comb.~29 (2009), 357-387

\bibitem{BoerGoer} J.~de Boer/J.~Goeree: Markov traces and $II_1$ factors
in conformal field theory, Commun.~Math.~Phys.~139 (1991), 267-304

\bibitem{DrinfQuasiHopf} V.~G.~Drinfel'd: Quasi-Hopf algebras, Leningrad Math.~J.~1 (1990), 1419-1457

\bibitem{Eholz1} W.~Eholzer: On the classification of modular fusion algebras, Commun.~Math.~Phys.~172 (1995), 623-659

\bibitem{EtVafa} P.~Etingof: On Vafa's theorem for tensor categories, 
Math.~Res.~Lett.~9 (2002), 651-657

\bibitem{EG4} P.~Etingof/S.~Gelaki: On the exponent of finite-dimensional Hopf algebras, Math.~Res.~Lett.~6 (1999), 131-140

\bibitem{EGQuasiHopf} P.~Etingof/S.~Gelaki: Finite dimensional quasi-Hopf algebras with radical of codimension~2, Math.~Res.~Lett.~11 (2004), 685-696

\bibitem{EtNikOst} P.~Etingof/D.~Nikshych/V.~Ostrik: On fusion categories,
Ann.~Math., II.~Ser.,~162 (2005), 581-642

\bibitem{FarbDennis} B.~Farb/R.~K.~Dennis: Noncommutative algebra, 
Grad.~Texts Math., Vol.~144, Springer, Berlin, 1993

\bibitem{Fine} B.~Fine: Algebraic theory of the Bianchi groups, Pure Appl.~Math., Vol.~129, Dekker, New York, 1989

\bibitem{FineRosen} B.~Fine/G.~Rosenberger: Number theory, Birkh\"auser, Basel, 2007

\bibitem{Frasch} H.~Frasch: Die Erzeugenden der Hauptkongruenzuntergruppen f\"ur Prim\-zahlstufen, Math.~Ann.~108 (1933), 229-252

\bibitem{GanModDat} T.~Gannon: Modular data: the algebraic combinatorics of conformal field theory, J.~Algebr.~Comb.~22 (2005), 211-250

\bibitem{GanMoonMonst} T.~Gannon: Moonshine beyond the monster, Camb.~Monogr.~Math.~Phys., Camb.~Univ.~Press, Cambridge, 2006

\bibitem{GaussSum} C.~F.~Gau\ss{}: Summatio quarumdam serierum singularium, Comm.\ Soc.\ reg.\ sci.\ Gottingensis rec.~1 (1811)

\bibitem{HausNillDoub} F.~Hau{\ss}er/F.~Nill: Doubles of quasi-quantum
groups, Commun.\ Math.\ Phys.~199 (1999), 547-589

\bibitem{Hu} S.~Hu: Lecture notes on Chern-Simons-Witten theory, World Scientific, Singapore, 2001

\bibitem{Hup1} B.~Huppert: Endliche Gruppen~I, Grundlehren Math.\ Wiss.,
Vol.~134, Springer, Berlin, 1979

\bibitem{Jac3} N.~Jacobson: Lectures in abstract algebra~III: Theory of fields and Galois theory, Grad.~Texts Math., Vol.~32, Springer, Berlin, 1976

\bibitem{K2} V.~G.~Kac: Infinite dimensional Lie algebras, 3rd~ed., Camb.~Univ.~Press, Cambridge, 1990

\bibitem{KashAnti} Y.~Kashina: On the order of the antipode of Hopf algebras in ${}_H^H YD$, Commun. Algebra 27 (1999), 1261-1273

\bibitem{KashPow} Y.~Kashina: A generalized power map for Hopf algebras. In: S.~Caene\-peel/ F.~van Oystaeyen (ed.): Hopf algebras and quantum groups, Lect.\ Notes Pure Appl.~Math., Vol.~209, Dekker, New York, 2000, 159-175  

\bibitem{YYY1} Y.~Kashina/Y.~Sommerh\"auser/Y.~Zhu: Self-dual modules of semisimple Hopf algebras, J.~Algebra~257 (2002), 88-96 

\bibitem{YYY2} Y.~Kashina/Y.~Sommerh\"auser/Y.~Zhu: On higher Frobenius-Schur indicators, Mem.~Am.~Math.~Soc., Vol.~181, No.~855, Am.~Math.~Soc., Providence, 2006

\bibitem{Kas} C.~Kassel: Quantum groups, Grad.~Texts Math., Vol.~155, Springer, Berlin, 1995

\bibitem{KasTur} C.~Kassel/V.~G.~Turaev: Braid groups, Grad.~Texts Math., Vol.~247, Springer, Berlin, 2008

\bibitem{KleinFricke} F.~Klein/R.~Fricke: Vorlesungen \"uber die Theorie der elliptischen Modulfunktionen, 1.~Band, Teubner, Stuttgart, 1890

\bibitem{Lan} E.~Landau: Elementary number theory, 2nd ed., Chelsea, New York, 1966

\bibitem{LR1} R.\ G.\ Larson/D.\ E.\ Radford: Semisimple cosemisimple Hopf algebras, Am.\ J.\ Math.\ 109 (1987), 187-195

\bibitem{LR2} R.~G.~Larson/D.~E.~Radford: Finite dimensional cosemisimple Hopf algebras in characteristic 0 are semisimple, J.~Algebra 117 (1988), 267-289

\bibitem{MasNg} G.~Mason/S.-H.~Ng: Central invariants and Frobenius-Schur indicators for semisimple quasi-Hopf algebras, Adv.~Math.~190 (2005), 161-195

\bibitem{Men} J.~Mennicke: On Ihara's modular group, Invent.\ Math.~4 (1967), 202-228

\bibitem{Mert} F.~Mertens: \"Uber die Gaussischen Summen, Berl. Ber.\ (1896), 217-219

\bibitem{M} S.~Montgomery: Hopf algebras and their actions on rings, 2nd 
revised printing, Reg.~Conf.~Ser.~Math., Vol.~82, Am.~Math.~Soc., Providence, 1997

\bibitem{MoSeib} G.~Moore/N.~Seiberg: Classical and quantum conformal field theory, Commun.~Math.~Phys.~123 (1989), 177-254

\bibitem{Mue1} M.~M\"uger: From subfactors to categories and topology I: Frobenius algebras in and Morita equivalence of tensor categories, J.~Pure Appl.~Algebra~180 (2003), 81-157 

\bibitem{Mue2} M.~M\"uger: From subfactors to categories and topology~II: The quantum double of tensor categories and subfactors, J.~Pure Appl.~Algebra 180 (2003), 159-219

\bibitem{Nag} T.~Nagell: Introduction to number theory, Chelsea, New York, 1964

\bibitem{NatExpFin} S.~Natale: On the exponent of tensor categories coming from finite groups, Israel J. Math. 162 (2007), 253-273

\bibitem{N} J.~Neukirch: Algebraische Zahlentheorie, Springer, Berlin, 1992

\bibitem{Newman} M.~Newman: Integral matrices, Pure Appl.~Math., Vol.~45, Academic Press, New York, 1972

\bibitem{NgSchauen2} S.-H.~Ng/P.~Schauenburg: Central invariants and higher Frobenius-Schur indicators for semisimple quasi-Hopf algebras, Trans.~Am.~Math.~Soc.~360 (2008), 1839-1860

\bibitem{NgSchauen3} S.-H.~Ng/P.~Schauenburg: Frobenius-Schur indicators and exponents of spherical categories, Adv.~Math.~211 (2007), 34-71

\bibitem{NgSchauen4} S.-H.~Ng/P.~Schauenburg: Congruence subgroups and generalized Frobe\-nius-Schur indicators, Preprint, arXiv:0806.2493v1 [math.QA], 2008

\bibitem{NR2} W.\ D.\ Nichols/M.\ B.\ Richmond: The Grothendieck algebra of a Hopf algebra I, Commun.~Algebra 26 (1998), 1081-1095 

\bibitem{RadElemNum} H.~Rademacher: Lectures on elementary number theory, Blaisdell, New York, 1964

\bibitem{RadfAntiQuasi} D.~E.~Radford: On the antipode of a quasitriangular Hopf algebra, J.~Algebra~151 (1992), 1-11

\bibitem{RadfKnotInv} D.~E.~Radford: On Kauffman's knot invariants arising from 
finite-dimen\-sional Hopf algebras. In: J.~Bergen/S.~Montgomery (ed.): Advances in Hopf algebras, Lect.~Notes Pure Appl.~Math., Vol.~158, Dekker, New York, 1994, 205-266

\bibitem{Rotman} J.~J.~Rotman: An introduction to homological algebra, Pure Appl.~Math., Vol.~85, Academic Press, New York, 1979

\bibitem{RowStongWang} E.~Rowell/R.~Stong/Z.~Wang: On classification of modular tensor categories, Preprint, arXiv:0712.1377v3 [math.QA], 2007

\bibitem{SchneiderFact} H.-J.~Schneider: Some properties of factorizable Hopf algebras, Proc.~Am.\ Math.~Soc.~129 (2001), 1891-1898

\bibitem{SchauenQuasiHopf} P.~Schauenburg: On the Frobenius-Schur indicators for quasi-Hopf algebras, J.~Algebra~282 (2004), 129-139

\bibitem{Serre2} J.~P.~Serre: Local fields, Grad.~Texts Math., Vol.~67, Springer, Berlin, 1979

\bibitem{So5} Y.~Sommerh\"auser: Ribbon transformations, integrals, and triangular decompositions, J.~Algebra~282 (2004), 423-489

\bibitem{SoZhu} Y.~Sommerh\"auser/Y.~Zhu: Hopf algebras and congruence subgroups, Preprint, arXiv:0710.0705v2 [math.RA], 2008

\bibitem{SoZhuRib} Y.~Sommerh\"auser/Y.~Zhu: On the notion of a ribbon quasi-Hopf algebra, Preprint, 2009

\bibitem{Schur} I.~Schur: \"Uber die Gaussschen Summen, G\"ott.~Nachr.\ (1921), 147-153 

\bibitem{Sw1} M.\ E.\ Sweedler: Integrals for Hopf algebras,
Ann.~Math., II.~Ser., 89 (1969), 323-335

\bibitem{Sw} M.\ E.\ Sweedler: Hopf algebras, W.\ A.\ Benjamin, New York, 1969

\bibitem{Tur} V.~G.~Turaev: Quantum invariants of knots and 3-manifolds, 
de Gruyter Stud.~Math., Vol.~18, de Gruyter, Berlin, 1994

\bibitem{Wash} L.~C.~Washington: Introduction to cyclotomic fields, Grad.~Texts Math., Vol.~83, Springer, Berlin, 1982

\bibitem{With} S.~J.~Witherspoon: The representation ring and the centre of a Hopf algebra, Can.~J.~Math.~51 (1999), 881-896 
\end{thebibliography}
\end{document}